# Interchange Rings
## Charles C. Edmunds

**Abstract**
An *interchange ring*, $(R,+,\bullet)$, is an abelian group with a second binary operation defined so that the *interchange law* $(x+y)\bullet(u+v) = (x\bullet u)+(y\bullet v)$ holds. An *interchange near ring* is the same structure based on a group which may not be abelian. It is shown that given any group, $G$, each interchange (near) ring based on that group is formed from a pair of endomorphisms of $G$ whose images commute, and that all interchange (near) rings based on $G$ can be characterized in this manner. To obtain an associative interchange ring, the endomorphisms must be commuting idempotents in the endomorphism semigroup of $G$. When $G$ is abelian we develop a group theoretic analogue of the simultaneous diagonalization of idempotent linear operators and show that pairs of endomorphisms which yield associative interchange rings can be diagonalized and then put into a canonical form. A best possible upper bound of $4^r$ can be given for the number of distinct isomorphism classes of associative interchange rings based on a finite abelian group $A$ which is a direct sum of $r$ cyclic groups of prime power order. If $A$ is direct sum of $r$ copies of the same cyclic group of prime power order, we show that there are exactly $\frac{1}{6}(r+1)(r+2)(r+3)$ distinct isomorphism classes of associative interchange rings based on $A$. Several examples are given and further comments are made about the general theory of interchange rings.

**1. Introduction.**
In [5] Kock introduced the notion of a double semigroup in his study of two-fold monoidal categories. A double semigroup is a nonempty set with two associative binary operations which satisfy the *interchange law*, (1.1). As early as 1962 Eckmann and Hilton [2] introduced the interchange law in their well-known argument demonstrating commutativity of the higher homotopy groups. Studies in double categories also make fruitful use of this law. Recently DeWolf [1] added to the study begun by Kock showing that all double inverse semigroups are commutative with both operations identical. In [3] the author discussed the construction of various double semigroups using commutation operations on a group. The purpose of this note is to study the implications of the interchange law in a universal algebraic context. Primarily, we will consider double magma whose first operation is a group and develop the obvious analogies with near rings and rings.

We will begin in the most general context and add further axioms as the theory develops. A *magma* $(M,*)$ is a pair consisting of a nonempty set and a binary operation on that set. A *double magma* $(M,*,\bullet)$ is a triple consisting of a nonempty set with two binary operations satisfying the *interchange law*,

$$(w*x)\bullet(y*z) = (w\bullet y)*(x\bullet z) \qquad (1.1)$$

for each $w,x,y,z \in M$. We call a double magma *commutative*, *associative*, or *unital* provided both magma, $(M,*)$ and $(M,\bullet)$, are, respectively, commutative, associative, or unital. An associative double magma is called a *double semigroup*. Expressed in this terminology, the Eckmann-Hilton argument proves that if a double magma is unital, then it is commutative, associative, the two identity elements for the operations coincide and, in fact, the two binary operations are the same. A double magma, $(M,*,*)$, with two identical operations is a double magma in name only; it coincides, essentially, with a magma, $(M,*)$, satisfying the *medial law*, $(w*x)*(y*z) = (w*y)*(x*z)$. If this magma is a semigroup, the medial law has a strong impact on its structure. For example if $(M,*)$ is finite or uniformly periodic, then it is commutative. In what follows we will try to avoid the situation where the two operations coincide; thus, as the Eckmann-Hilton theorem implies, we must not allow a double magma to be unitary. Henceforth we will allow at most one of the operations to have an identity element. With this point in mind, we define a double magma to be *proper* if its operations do not coincide and *improper* otherwise.

The specific focus of this study is an investigation of the analogy between a double magma $(M,*,\bullet)$ and a ring $(R,+,\cdot)$. A ring is an additive abelian group with an associative multiplication and distributive laws serving to connect the two operations algebraically. We will not consider rings with identity since, the analogous double magma is then forced to be improper. A double magma has no algebraic structure on its two binary operations, but, as with the distributive laws in a ring, the interchange law connects the two operations. A *near ring* is a generalization of a ring in which the additive structure may not be abelian. Also only one of the distributive laws, usually the right distributive law, is imposed for a near ring. We begin our study with an arbitrary group $(G,+)$ written additively using the usual conventions that the identity element is written 0 and the inverse of $x \in G$ is written $-x$. We will form a double magma $(G,+,\bullet)$ where the magma $(G,\bullet)$ has no algebraic structure beyond the closure of its binary operation on $G$. Thus, with the distributive laws replaced by the interchange law, we have something like a near ring, but with the associative law not necessarily holding on the "multiplication".

**Definition:** An *interchange near ring* is a triple, $(G,+,\bullet)$, with $(G,+)$ a group and $(G,+,\bullet)$ a double magma.

If the additive structure is based on an abelian group, we specialize the definition.

**Definition:** An *interchange ring* is an interchange near ring, $(G,+,\bullet)$, for which $(G,+)$ is abelian.

Interchange rings may seem like exotic objects, but it is interesting to note that the well-know structures $(\mathbb{Z},+,-)$ and $(\mathbb{Q}\setminus\{0\},\times,\div)$ are interchange rings. In fact, if $(A,+)$ is any abelian group, $(A,+,-)$ is an interchange ring. It is a well-known algebraic problem to determine which abelian groups admit a non-zero ring structure. For example the Prufer group, or quasicyclic group $\mathbb{Z}_{p^\infty}$, can be seen to admit only the zero ring structure. In Section 3 we show that every nontrivial group, $G$, admits at least three nonzero interchange (near) ring structures and there is a clearly defined process for constructing each such structure from the group. The construction starts with a pair, $(\varepsilon,\eta)$, of endomorphisms of $G$ with commuting images and then defines a product by $x \bullet y = \varepsilon(x) + \eta(y)$. It is shown, Theorem 3.2, that each product constructed in this way yields an interchange near ring, and, Theorem 3.3, that each interchange near ring based on $G$ arises in this manner. It is then observed, Theorem 3.5, that the equivalence relation of similarity of pairs, i.e. the existence of $\alpha \in Aut(G,+)$ for which $\alpha^{-1}(\varepsilon_1,\eta_1)\alpha = (\varepsilon_2,\eta_2)$, corresponds exactly to isomorphism of the associated interchange (near) rings.

In Section 4 we show, Theorem 4.3, that an interchange (near) ring is associative if and only if it is generated by a pair of commuting, idempotent endomorphisms of its additive group. In Section 5 we restrict our attention to interchange rings based on finite abelian groups. We apply techniques familiar in operator theory to characterize all finite associative interchange rings by showing how to 'diagonalize' pairs of their generating endomorphisms, Theorem 5.9. We obtain, Theorem 5.10, a best possible upper bound of $4^r$ for the number of distinct isomorphism classes of associative interchange rings based on a finite abelian group, $A$, which is a direct sum of $r$ cyclic groups of prime power order. In Section 6 we further restrict $A$ to a direct sum of $r$ isomorphic copies of a cyclic group of prime power order and develop a unique canonical form, Corollary 6.8, for pairs of endomorphism generating the interchange rings based on $A$. Using this result, we obtain, Theorem 6.9, the exact number, $\frac{1}{6}(r+1)(r+2)(r+3)$, of distinct isomorphism classes of associative interchange rings based on $A$. In Section 7 we make some comments about the extent to which interchange ring theory might be developed by analogy with standard ring theory. One would hope, one day, to see a structure theory of interchange rings.

**2. Preliminaries.**
We begin with some fundamental statements about interchange near rings. Note the algebraic importance of the zero element.

**Lemma 2.1:** If $(G,+,\bullet)$ is an interchange near ring, then the following statements are true.
(i) Zero is idempotent multiplicatively; $0 \bullet 0 = 0$.
(ii) Zero distributes multiplicatively over addition; for each $x, y \in G$,
$$0 \bullet (x+y) = (0 \bullet x) + (0 \bullet y) \text{ and } (x+y) \bullet 0 = (x \bullet 0) + (y \bullet 0).$$
(iii) For each $x, y \in G$, $x \bullet y = x \bullet 0 + 0 \bullet y$.
(iv) For each $x, y \in G$, $(-x) \bullet (-y) = -(x \bullet y)$.

**Proof: (i)** We use the identity property of 0 and then apply the interchange law to obtain,
$$0 \bullet 0 = (0+0) \bullet (0+0) = (0 \bullet 0) + (0 \bullet 0).$$
Adding the inverse of $0 \bullet 0$ to both sides, we obtain our result.
**(ii)** $0 \bullet (x+y) = (0+0) \bullet (x+y) = (0 \bullet x) + (0 \bullet y)$. The second equality follows by the interchange law. The right-hand law follows similarly.
**(iii)** $x \bullet y = (x+0) \bullet (0+y) = x \bullet 0 + 0 \bullet y$.
**(iv)** $(x \bullet y) + (-x \bullet -y) = (x+(-x)) \bullet (y+(-y)) = 0 \bullet 0 = 0$. The result follows. ∎

**Lemma 2.2:** If $(G,+,\bullet)$ is an interchange near ring then $(G,\bullet)$ is associative if and only if, for each $x \in G$,
(i) $(x \bullet 0) \bullet 0 = x \bullet 0$, (ii) $(0 \bullet x) \bullet 0 = 0 \bullet (x \bullet 0)$, and (iii) $0 \bullet (0 \bullet x) = 0 \bullet x$.

**Proof:** Note first that if $(G,\bullet)$ is associative then (ii) follows immediately. Conditions (i) and (ii) follows from associativity and Lemma 2.1(i). Conversely, assuming (i), (ii), and (iii) hold and applying Lemma 2.1,
$$(x \bullet y) \bullet z = (x \bullet 0 + 0 \bullet y) \bullet 0 + 0 \bullet z = (x \bullet 0) \bullet 0 + (0 \bullet y) \bullet 0 + 0 \bullet z = x \bullet 0 + 0 \bullet (y \bullet 0) + 0 \bullet (0 \bullet z)$$
$$= x \bullet 0 + 0 \bullet ((y \bullet 0) + (0 \bullet z)) = x \bullet (y \bullet z). \blacksquare$$

In rings the additive identity is an annihilator (that is $x \cdot 0 = 0 \cdot x = 0$, for each $x$), however the existence of an annihilator in an interchange near ring destroys its structure. We define an element $a \in G$ to be an *annihilator* in an interchange near ring, $(G,+,\bullet)$, if $a \bullet x = x \bullet a = a$ for every $x \in G$. A *zero semigroup* is a semigroup $(S,\cdot)$ satisfying the identity $w \cdot x = y \cdot z$. This implies that all entries in the Cayley table are equal; thus zero semigroups are considered trivial for most purposes.

**Proposition 2.3:** If $(G,+,\bullet)$ is an interchange near ring containing an annihilator, then $(G,\bullet)$ is a zero semigroup.

**Proof:** Suppose $a$ is an annihilator. Thus, for each $x \in G$, we have $a = a \bullet x = a \bullet 0 + 0 \bullet x = a + 0 \bullet x$. Therefore $0 \bullet x = 0$, for each $x \in G$. Similarly we can show that $x \bullet 0 = 0$, for each $x \in G$. Thus, by Lemma 2.1 (i) and (iii), we have $x \bullet y = x \bullet 0 + 0 \bullet y = 0 + 0 = 0$, for each $x, y \in G$. Hence $(G, \bullet)$ is a zero semigroup. ∎

## 3. Constructing Interchange Near Rings from Groups.

In this section we will show that any nontrivial group admits a number of distinct interchange (near) ring structures. For $(G, +)$ a group, we denote the sets of endomorphisms and automorphisms of $G$ by $End(G, +)$ and $Aut(G, +)$, respectively.

**Definition:** If $(G, +)$ is a group with $\varepsilon, \eta \in End(G, +)$, we say the pair $(\varepsilon, \eta)$ is *image-commuting* when $\varepsilon(x) + \eta(y) = \eta(y) + \varepsilon(x)$, for every $x, y \in G$.

If $(G, +)$ is abelian, then each pair of endomorphisms is image commutating. If $(G, +)$ is nonabelian and the image of $G$ under either map is contained in the centre of $G$, then the maps are image-commuting. We give an example to illustrate that image-commuting endomorphisms need not have abelian images.

**Example 3.1:** We present the symmetric group on three letters additively as
$$S_3 = \langle a, b; 3a = 0, 2b = 0, b + a = 2a + b \rangle$$
Each element of $S_3$ can be written uniquely as $ia + jb$ with $i \in \mathbb{Z}_3$ and $j \in \mathbb{Z}_2$. Let $(G, +)$ be the direct product of two copies of $S_3$, and define two endomorphisms of $(G, +)$ as follows: for each $(ia + jb, ka + lb) \in G$, let $\varepsilon(ia + jb, ka + lb) = (0, ia + jb)$ and let $\eta(ia + jb, ka + lb) = (ka + lb, 0)$. Note that the images of both $\varepsilon$ and $\eta$ are isomorphic to $S_3$, and are, hence, nonabelian. However, since the images of these endomorphisms are in different direct factors, they commute; thus $(\varepsilon, \eta)$ is an image-commuting pair.

We are now in a position to state the main results of this section.

**Theorem 3.2:** If $(G, +)$ is a group, $(\varepsilon, \eta)$ is an image-commuting pair of endomorphisms of $(G, +)$, and a binary operation $\bullet$ is defined on $G$ by $x \bullet y = \varepsilon(x) + \eta(y)$, for each $x, y \in G$, then $(G, +, \bullet)$ in an interchange near ring.

And conversely, we have the following.

**Theorem 3.3:** If $(G, +, \bullet)$ is an interchange near ring, then the pair $(\varepsilon, \eta)$ of mappings from $G$ to $G$ defined by $\varepsilon(x) = x \bullet 0$ and $\eta(x) = 0 \bullet x$, for each $x \in G$, is a pair of image-commutating endomorphisms of $G$.

We will denote the set of all image-commuting pairs of endomorphisms of $(G, +)$ by $ICE(G, +)$, and the set of all interchange near rings based on $(G, +)$ by $INR(G, +)$. Given a pair $(\varepsilon, \eta) \in ICE(G, +)$, we denote the binary operation constructed from this pair in Theorem 3.2 by $\bullet_{(\varepsilon, \eta)}$. Given an interchange near ring, $(G, +, \bullet)$, we denote the pair of image-commuting mappings defined in Theorem 3.3 by $(\varepsilon_\bullet, \eta_\bullet)$.

**Theorem 3.4:** The mapping $\Psi : ICE(G) \to INR(G)$ defined by $\Psi(\varepsilon, \eta) = (G, +, \bullet_{(\varepsilon,\eta)})$ is a bijection.

We will leave the task of determining structural parallels between these two classes to category theorists.

**Proof of Theorem 3.2:** It is sufficient to show that the interchange law holds. Letting $w, x, y, z \in G$, we have $(w \bullet x) + (y \bullet z) = \varepsilon(w) + \eta(x) + \varepsilon(y) + \eta(z)$. Since $\varepsilon$ and $\eta$ are image-commuting, we can interchange the two middle terms, and obtain

$$(w \bullet x) + (y \bullet z) = \varepsilon(w) + \varepsilon(y) + \eta(x) + \eta(z) = \varepsilon(w + y) + \eta(x + z) = (w + y) \bullet (x + z). \blacksquare$$

**Proof of Theorem 3.3:** Clearly $\varepsilon(x) = x \bullet 0$ maps from $G$ into $G$. To show that $\varepsilon$ is a homomorphism, we let $x, y \in G$ and consider $\varepsilon(x + y)$ and $\varepsilon(-x)$. Applying Lemma 2.1 we have, $\varepsilon(x + y) = (x + y) \bullet 0 = x \bullet 0 + y \bullet 0 = \varepsilon(x) + \varepsilon(y)$, and $\varepsilon(-x) = (-x) \bullet 0 = -(x \bullet 0) = -\varepsilon(x)$. Thus $\varepsilon$ is a homomorphism of $G$. The proof that $\eta$ is a homomorphism of $G$ is similar. Lastly, we must see that the pair is image-commuting. Let $x, y \in G$ and consider the sum, $\varepsilon(x) + \eta(y)$. By Lemma 2.1 and the interchange law we have,

$$\varepsilon(x) + \eta(y) = x \bullet 0 + 0 \bullet y = x \bullet y = (0 + x) \bullet (y + 0) = (0 \bullet y) + (x \bullet 0) = \eta(y) + \varepsilon(x). \blacksquare$$

**Proof of Theorem 3.4:** First we will show that $\Psi$ maps $ICE(G)$ into $INR(G)$. It suffices to show that the image of any $(\varepsilon, \eta) \in ICE(G, +)$ under $\Psi$, namely $(G, +, \bullet_{(\varepsilon,\eta)})$, satisfies the interchange law. Supposing that $w, x, y, z \in G$, we have,

$$(w \bullet_{(\varepsilon,\eta)} x) + (y \bullet_{(\varepsilon,\eta)} z) = (\varepsilon(w) + \eta(x)) + (\varepsilon(y) + \eta(z)).$$

Since $(\varepsilon, \eta)$ is an image-commuting pair, we can rearrange the second sum as

$$(\varepsilon(w) + \eta(x)) + (\varepsilon(y) + \eta(z)) = (\varepsilon(w) + \varepsilon(y)) + (\eta(x) + \eta(z)).$$

Since $\varepsilon$ and $\eta$ are endomorphisms of G, we can apply the definition of $\bullet_{(\varepsilon,\eta)}$ in reverse, to conclude that,

$$(\varepsilon(w) + \varepsilon(y)) + (\eta(x) + \eta(z)) = \varepsilon(w + y) + \eta(x + z) = (w + y) \bullet_{(\varepsilon,\eta)} (x + z).$$

To prove $\Psi$ is injective, we suppose $(\varepsilon_1, \eta_1)$ and $(\varepsilon_2, \eta_2)$ are image-commuting pairs for which $\Psi(\varepsilon_1, \eta_1) = \Psi(\varepsilon_2, \eta_2)$ that is, $(G, +, \bullet_{(\varepsilon_1,\eta_1)}) = (G, +, \bullet_{(\varepsilon_2,\eta_2)})$, and therefore $\bullet_{(\varepsilon_1,\eta_1)} = \bullet_{(\varepsilon_2,\eta_2)}$. To keep our calculations more readable, we will use $\bullet_i$ as shorthand notation for $\bullet_{(\varepsilon_i,\eta_i)}$. Since each $\eta_i$ is a homomorphism of $G$, it follows that $\eta_i(0) = 0$. Thus, for any $x \in G$,

$$\varepsilon_1(x) = \varepsilon_1(x) + \eta_1(0) = x \bullet_1 0 = x \bullet_2 0 = \varepsilon_2(x) + \eta_2(0) = \varepsilon_2(x).$$

Thus $\varepsilon_1 = \varepsilon_2$ and a similar calculation shows that $\eta_1 = \eta_2$. Therefore $(\varepsilon_1, \eta_1) = (\varepsilon_2, \eta_2)$ and, thus, $\Psi$ is injective. To see that $\Psi$ is surjective, we let $(G, +, \bullet)$ be any element of $INR(G, +)$. Theorem 3.3 shows how to construct the pair, $(\varepsilon_\bullet, \eta_\bullet)$, of image-commuting endomorphisms of $(G, +)$ from $(G, +, \bullet)$. Applying $\Psi$ to this pair, we obtain the interchange near ring $\Psi(\varepsilon_\bullet, \eta_\bullet)$ based on $(G, +)$. For now, let us denote the second operation of the interchange near ring thus produced by $\odot$, writing $\Psi(\varepsilon_\bullet, \eta_\bullet) = (G, +, \odot)$. Note that $\odot$ is constructed from the pair $(\varepsilon_\bullet, \eta_\bullet)$ as in Theorem 3.2. That is, for each $x, y \in G$,

$$x \odot y = \varepsilon_\bullet(x) + \eta_\bullet(y) = x \bullet 0 + 0 \bullet y = x \bullet y.$$

Therefore $\odot = \bullet$, and as a result, $(G, +, \odot) = (G, +, \bullet)$. It follows that $\Psi(\varepsilon_\bullet, \eta_\bullet) = (G, +, \bullet)$, and hence, that $\Psi$ is surjective. $\blacksquare$

At this point we have shown that distinct pairs of image-commuting endomorphisms yield unequal interchange near rings, but they could be isomorphic. If $(G_1,+_1,\bullet_1)$ and $(G_2,+_2,\bullet_2)$ are interchange near rings, a mapping
$$\varphi:(G_1,+_1,\bullet_1) \to (G_2,+_2,\bullet_2),$$
is an *interchange near ring homomorphism* if $\varphi:(G_1,+_1) \to (G_2,+_2)$ is a group homomorphism and, for each $x,y \in G_1$, $\varphi(x \bullet_1 y) = \varphi(x) \bullet_2 \varphi(y)$. This homomorphism is an *interchange near ring isomophism* exactly when it is a bijection. We will denote isomorphism of interchange near rings by $\cong$. It is routine to check that both relations defined below are equivalence relations.

**Definition:** If $\varepsilon,\eta \in End(G,+)$ and $\alpha \in Aut(G,+)$ so that $\alpha^{-1}\varepsilon\alpha = \eta$, we say $\varepsilon$ is similar to $\eta$ (under $\alpha$) and write $\varepsilon \sim \eta(\alpha)$ or, more simply, $\varepsilon \sim \eta$. For pairs of endomorphisms of $G$, $(\varepsilon_1,\eta_1)$ is similar to $(\varepsilon_2,\eta_2)$, denoted $(\varepsilon_1,\eta_1) \sim (\varepsilon_2,\eta_2)(\alpha)$ or $(\varepsilon_1,\eta_1) \sim (\varepsilon_2,\eta_2)$ whenever there is an automorphism $\alpha \in Aut(G,+)$ such that $\alpha^{-1}(\varepsilon_1,\eta_1)\alpha = (\alpha^{-1}\varepsilon_1\alpha, \alpha^{-1}\eta_1\alpha) = (\varepsilon_2,\eta_2)$.

**Theorem 3.5:** If $(G,+)$ is a group and $(G,+,*)$ and $(G,+,\odot)$ are interchange near rings based on $(G,+)$, then $(G,+,*) \cong (G,+,\odot)$ if and only $(\varepsilon_*,\eta_*) \sim (\varepsilon_\odot,\eta_\odot)$.

**Proof:** ($\Rightarrow$) Suppose there is an isomorphism $\varphi:(G,+,*) \to (G,+,\odot)$. First consider $\varepsilon_*$ applied to any $x \in G$, $\varepsilon_*(x) = x*0$. Applying the isomorphism $\varphi$ we obtain,
$$\varphi(\varepsilon_*(x)) = \varphi(x*0) = \varphi(x) \odot \varphi(0) = \varphi(x) \odot 0 = \varepsilon_\odot(\varphi(x)).$$
Therefore, $(\varphi\varepsilon_* - \varepsilon_\odot\varphi)(x) = 0$, for each $x \in G$. It follows that $\varphi\varepsilon_* = \varepsilon_\odot\varphi$, and by a similar argument that $\varphi\eta_* = \eta_\odot\varphi$. Thus we have $(\varepsilon_*,\eta_*) \sim (\varepsilon_\odot,\eta_\odot)$.

($\Leftarrow$) Now suppose that $(\varepsilon_*,\eta_*) \sim (\varepsilon_\odot,\eta_\odot)$. Thus there is an automorphism $\phi$ of $(G,+)$ such that $\phi\varepsilon_* = \varepsilon_\odot\phi$ and $\phi\eta_* = \eta_\odot\phi$. We claim that $\phi$ is an isomorphism mapping $(G,+,*)$ onto $(G,+,\odot)$. Since $\phi$ is an automorphism of $(G,+)$ we know that it is a bijection from $G$ to $G$ and that it is a group homomorphism with respect to addition. To prove that $(G,+,*)$ and $(G,+,\odot)$ are isomorphic under $\phi$, it remains to show that, for each $x,y \in G$, $\phi(x*y) = \phi(x) \odot \phi(y)$.
$$\phi(x*y) = \phi(x*0 + 0*y) = \phi(x*0) + \phi(0*y) = \phi\varepsilon_*(x) + \phi\eta_*(y) = \varepsilon_\odot\phi(x) + \eta_\odot\phi(y) = \phi(x) \odot \phi(y). \blacksquare$$

**Theorem 3.6:** If $(G,+)$ is a non-trivial group, there are at least three distinct isomorphism classes of interchange near rings based on $G$.

**Proof:** Let $\iota,\zeta$ be the identity map and the zero map from $G$ to $G$, respectively. Note that $\zeta$ is image-commuting with all endomorphisms of $G$. Since both of these mappings commute with every automorphism of $G$, their similarity classes are singletons. Thus, by Theorems 3.2 and 3.5, the pairs $(\zeta,\zeta),(\zeta,\iota)$, and $(\iota,\zeta)$, generate three non-isomorphic interchange near rings based on $G$. $\blacksquare$

Note that $(\zeta,\zeta)$ generates the product $x \bullet y = 0$, yielding the zero interchange ring, the pair $(\zeta,\iota)$ generates the product $x \bullet y = x$, the right zero semigroup, and $(\iota,\zeta)$ generates the product $x \bullet y = y$, the left zero semigroup. We will refer to these three isomorphism classes of interchange (near) rings as the *essential* interchange rings based on $G$. We will find these relatively uninteresting and refer to interchange (near) rings based on G which are not of these four types as *inessential*. The question, in ring theory, of whether or not there exists a non-zero ring based on a particular abelian group is parallel to the question, in interchange (near) ring theory, of whether or not there is an inessential interchange (near) rings based on a group.

With this much structure involved, the discussion calls for an example of how our theorems can be applied in some specific cases.

**Example 3.7: Interchange near rings based on $S_3$.** For economy, we will rename the elements of $(S_3,+)$ as given in Example 3.1, $0, a, 2a, b, a+b, 2a+b$, as 0,1,2,3,4,5, respectively. We denote an endomorphism $\varepsilon$ mapping $0 \mapsto 0, 1 \mapsto v, 2 \mapsto w, 3 \mapsto x, 4 \mapsto y,$ and $5 \mapsto z$ by $\varepsilon = (0vwxyz)$. There are six automorphisms of $(S_3,+)$,

$$\alpha_0 = (012345), \alpha_1 = (012453), \alpha_2 = (012534), \alpha_3 = (021354), \alpha_4 = (021435), \text{ and } \alpha_5 = (021543),$$

and four proper endomorphisms,

$$\varepsilon_0 = (000000), \varepsilon_1 = (000333), \varepsilon_2 = (000444), \text{ and } \varepsilon_3 = (000555).$$

The images of $S_3$ under these maps are

$$\alpha_i(S_3) = S_3, \varepsilon_0(S_3) = \{0\}, \varepsilon_1(S_3) = \{0,3\}, \varepsilon_2(S_3) = \{0,4\}, \text{ and } \varepsilon_3(S_3) = \{0,5\}.$$

To form an interchange near ring, the pairs we select must be image-commuting, thus the only candidates are pairs of the form $(\varepsilon_0, \varepsilon)$ and $(\varepsilon, \varepsilon_0)$ for $\varepsilon \in End(S_3,+)$, and the pairs $(\varepsilon_1, \varepsilon_1), (\varepsilon_2, \varepsilon_2)$, and $(\varepsilon_3, \varepsilon_3)$. Since $\alpha_1^{-1}\varepsilon_1\alpha_1 = \varepsilon_2$ and $\alpha_3^{-1}\varepsilon_2\alpha_3 = \varepsilon_3$, we have $(\varepsilon_1, \varepsilon_1) \sim (\varepsilon_2, \varepsilon_2) \sim (\varepsilon_3, \varepsilon_3)$; therefore, these last three pairs will yield isomorphic interchange rings. We express $End(S_3)$ as a disjoint union of similarity classes as:

$$End(S_3) = \{\alpha_0\} \dot{\cup} \{\alpha_1, \alpha_2\} \dot{\cup} \{\alpha_3, \alpha_4, \alpha_5\} \dot{\cup} \{\varepsilon_0\} \dot{\cup} \{\varepsilon_1, \varepsilon_2, \varepsilon_3\}.$$

Thus the non-similar, image-commuting pairs are $(\varepsilon_1, \varepsilon_1)$, $(\varepsilon_0, \varepsilon)$ where $\varepsilon \in \{\alpha_0, \alpha_1, \alpha_3, \varepsilon_0, \varepsilon_1\}$, and $(\varepsilon, \varepsilon_0)$ where $\varepsilon \in \{\alpha_0, \alpha_1, \alpha_3, \varepsilon_1\}$. Hence there are ten distinct isomorphism classes of interchange near rings based on $S_3$; the three essential interchange near rings generated by the pairs $(\varepsilon_0, \varepsilon_0), (\varepsilon_0, \alpha_0)$, and $(\alpha_0, \varepsilon_0)$ plus six additional inessential interchange near rings.

**4. Commutativity, Idempotence, and Associativity.**

**Proposition 4.1:** Given an interchange near ring $(G, +, \bullet)$ constructed from a pair, $(\varepsilon, \eta)$ of image-commuting endomorphisms of $(G,+)$, $(G, \bullet)$ is commutative if and only if $\varepsilon = \eta$.
**Proof:** Note first that if $(G, \bullet)$ is commutative, we have $\varepsilon(x) + \eta(y) = x \bullet y = y \bullet x = \varepsilon(y) + \eta(x)$. Our conclusion follows letting $y = 0$. Conversely, if we suppose that $\varepsilon = \eta$, then, recalling that $\varepsilon$ and $\eta$ are image-commuting, we have $x \bullet y = \varepsilon(x) + \eta(y) = \eta(y) + \varepsilon(x) = \varepsilon(y) + \eta(x) = y \bullet x$. ∎

Note that when $(G, +, \bullet)$ is commutative, to be image-commuting $\varepsilon(G)(= \eta(G))$ must be abelian. We denote the identity mapping of $End(G,+)$ by $\iota$.

**Proposition 4.2:** Given an interchange near ring $(G, +, \bullet)$ constructed from a pair, $(\varepsilon, \eta)$ of image-commuting endomorphisms of $(G,+)$, $(G, \bullet)$ is idempotent if and only if $\varepsilon + \eta = \iota$.

**Proof:** Note that $x = x \bullet x = \varepsilon(x) + \eta(x) = (\varepsilon + \eta)(x).$ It follows that idempotence of $(G, \bullet)$ is equivalent to $\varepsilon + \eta = \iota.$ ∎

The following theorem characterizes those interchange near rings which have an associative multiplication.

**Theorem 4.3:** Given an interchange near ring $(G, +, \bullet)$ constructed from a pair, $(\varepsilon, \eta)$ of image-commuting endomorphisms of $(G, +)$, $(G, \bullet)$ is associative if and only if $\varepsilon$ and $\eta$ are commuting idempotents in the semigroup $End(G, +)$.

**Proof:** By Lemma 2.2, $(G, \bullet)$ is associative if and only if, for each $x \in G$, (i) $(x \bullet 0) \bullet 0 = x \bullet 0$, (ii) $(0 \bullet x) \bullet 0 = 0 \bullet (x \bullet 0)$, and (iii) $0 \bullet (0 \bullet x) = 0 \bullet x$. Since $\varepsilon(x) = x \bullet 0$, we have $\varepsilon^2(x) = \varepsilon(\varepsilon(x)) = \varepsilon(x) \bullet 0 = (x \bullet 0) \bullet 0$. Therefore condition (i) is equivalent to the idempotence of $\varepsilon$. Similarly condition (iii) is equivalent to the idempotence of $\eta$. Finally, note that $(0 \bullet x) \bullet 0 = \varepsilon(0 \bullet x) = \varepsilon(\varepsilon(0) + \eta(x)) = \varepsilon\eta(x)$, while $0 \bullet (x \bullet 0) = \varepsilon(0) + \eta(x \bullet 0) = \eta(\varepsilon(x) + \eta(0)) = \eta\varepsilon(x)$. Thus condition (ii) is equivalent to the commutativity of $\varepsilon$ and $\eta$. ∎

**Example 4.4: Associative interchange near rings of order six:** The two groups of order six are the cyclic group of order 6, which we will present as $C_6 = \langle a; 6a = 0 \rangle$, and the symmetric group on three letters, $S_3$, which was discussed in Example 3.7. We will find 16 isomorphism classes of associative interchange rings based on $C_6$ and six isomorphism classes of associative interchange near rings with additive structure $(S_3, +)$. Thus we have a total of 22 distinct isomorphism classes of associative interchange near rings of order 6. Any endomorphism of $C_6$ can be written as $\varepsilon_i : x \mapsto ix$ for $i \in \mathbb{Z}_6$. The mappings $\varepsilon_1$ and $\varepsilon_5$ are automorphisms and the other mappings are proper endomorphisms. Since $C_6$ is abelian, all pairs of endomorphisms are image-commuting. It is easy to see that $(End(C_6, +), \circ)$ is a commutative semigroup with idempotents $E = \{\varepsilon_0, \varepsilon_1, \varepsilon_3, \varepsilon_4\}$. Thus for each $\alpha \in Aut(C_6, +)$ and for each $\varepsilon \in End(C_6, +)$, $\alpha^{-1}\varepsilon\alpha = \alpha$. It follows that each similarity class in $End(C_6, +)$ is a singleton, and, hence, each of the sixteen pairs of elements from $E$ generates a distinct isomorphism class of associative interchange rings. By Proposition 4.1, we see that the four pairs $(\varepsilon_i, \varepsilon_i)$ $(1 \le i \le 4)$ yield commutative multiplications. And by Proposition 4.2, the pairs $(\varepsilon_0, \varepsilon_1), (\varepsilon_1, \varepsilon_0), (\varepsilon_3, \varepsilon_4),$ and $(\varepsilon_4, \varepsilon_3)$ yield idempotent multiplications. Turning to $S_3$, we see from Example 3.7 that among the ten pairs which yield interchange near rings, there are six which are pairs of commuting idempotents:

$$(\varepsilon_0, \varepsilon_0), (\varepsilon_0, \alpha_0), (\alpha_0, \varepsilon_0), (\varepsilon_0, \varepsilon_1), (\varepsilon_1, \varepsilon_0), (\varepsilon_1, \varepsilon_1).$$

These generate six distinct isomorphism classes. By Proposition 4.1, we see that only the pairs $(\varepsilon_0, \varepsilon_0)$ and $(\varepsilon_1, \varepsilon_1)$ yield commutative multiplications. And by Proposition 4.2, the pairs $(\varepsilon_0, \alpha_0)$ and $(\alpha_0, \varepsilon_0)$ yield idempotent multiplications.

**5. Finite associative interchange rings.**

In this section we refer the reader to Chapter 8 of McCoy [6] for reference to results on abelian groups and to Chapter 6 of Hoffman and Kunze [4] for reference to linear algebra. We begin by proving some facts about endomorphisms of abelian groups which are analogous to familiar results in linear algebra, and facilitated by the fact that an abelian group is a $\mathbb{Z}$-module. The goal of this section, Theorem 5.9, is to apply the theory thus developed to give a characterization of all finite associative interchange rings. Given a finite abelian group $(A,+)$, we will show that an associative interchange ring based on $A$ is defined by a pair of 'diagonalizable' endomorphisms, and that, when all such pairs are made diagonal, up to similarity, they yield all isomorphism class of associative interchange ring based on $A$. This will be accomplished by characterizing the commuting, idempotent elements of $End(A,+)$ and then appealing to Theorem 4.3. We will decompose $A$ as a direct sum in a way that allows us to take advantage of an analogy with linear algebra. Luckily, it turns out, the theorem of linear algebra which we will need, that commuting, idempotent linear operators are simultaneously diagonalizable, survives in the context of abelian groups, Theorem 5.9. This diagonalization will allows us to give, Theorem 5.10, an upper bound, which is tight in a sense, on the number of associative interchange rings based on a finite abelian group $A$.

For the remainder of this section $(A,+)$ will be a finite additive abelian group with identity element 0. By the fundamental structure theorem, we know that $A$ can be decomposed as a direct sum of prime power order cyclic groups. We will refer to any such decomposition as a *ppc-decomposition*. Furthermore if $A$ is decomposed in this manner as both $A = A_1 \oplus A_2 \oplus \cdots \oplus A_s$ and $A = B_1 \oplus B_2 \oplus \cdots \oplus B_t$, then $s = t$ and there is bijection between the sets of summands so that each $A_i$ is isomorphic to the corresponding $B_j$. Since the number of these prime power order cyclic summands is invariant for $A$, we call this number the *prime power cyclic rank* of $A$, denoted *ppc-rank(A)*. The true invariant here is the multiset of prime power order, cyclic summands which form $A$ as a direct sum. Routine use of the fundamental structure theory for finite abelian groups shows that (i) ppc-rank is invariant under isomorphism, and (ii) if $A$ is a finite abelian group and $S \leq A$, then ppc-rank$(S) \leq$ ppc-rank$(A)$.

For each $i$ ($1 \leq i \leq r$) let $a_i \in A$ so that $A = \langle a_1 \rangle \oplus \langle a_2 \rangle \oplus \cdots \oplus \langle a_r \rangle$ with each summand a cyclic group of prime power order, and let $\bar{a} = \{a_1, \ldots, a_r\}$. Any $a \in A$ can be written uniquely as $a = c_1 a_1 + \cdots + c_r a_r$ with each $c_i \in \mathbb{Z}_{|a_i|}$, where $|a_i|$ denotes the order of $a_i$. It is reasonable to think of $\bar{a}$ as a basis for $A$, analogous to the standard concept for finite dimensional vector spaces. We will call any such set $\bar{a}$ a *prime power cyclic basis* of $A$, and abbreviate this term as *ppc-basis*. The trivial group has no ppc-basis and is said to have *ppc-rank* zero. We comment that the ppc-rank is not the usual rank of the abelian group. In fact a cyclic group of order 6 has rank 1 and ppc-rank 2.

If $\bar{x} = \{x_1, \ldots, x_r\}$ is a ppc-basis for $A$ and $\phi^* : \bar{x} \to A$ is a mapping for which each $|\phi^*(x_i)|$ divides $|x_i|$, then we can extend $\phi^*$ uniquely to an endomorphism $\phi \in End(A,+)$ linearly; that is, if $a \in A$ is written $a = c_1 x_1 + \cdots c_r x_r$ with each $c_i \in \mathbb{Z}_{|x_i|}$, we let $\phi(a) = c_1 \phi^*(x_1) + \cdots c_r \phi^*(x_r)$. The fact that this is a homomorphism can be seen by appealing to the presentation of $A$ implied by the ppc-basis $\bar{x}$. If we denote the order of each $x_i$ in $A$ as $n_i$, then $A$ can be presented as $A = \langle G; R \rangle$ where the generating set is $G = \bar{x}$ and the set of relations is $R = \{n_i x_i = 0 : 1 \leq i \leq r\} \cup \{x_i + x_j = x_j + x_i : 1 \leq i, j \leq r\}$. Any mapping of $G$ into $A$ extends to a unique homomorphism exactly when the mapping preserves these relations. Since we are mapping into an abelian group, the second set of relations holds under $\phi$. And the first set holds, that is $\phi(n_i x_i) = n_i \phi^*(x_i) = 0$, since we have made the hypothesis that $|\phi^*(x_i)|$ divides $|x_i|$.

For the remainder of Sections 5 and 6 we will let $A$ be a finite abelian group of ppc-rank$(A) = r$ and fix a particular pcc-basis $\bar{e} = \{e_1, \ldots, e_r\}$ for $A$, which we will refer to as the *standard pcc-basis* for $A$. We say an endomorphism $\delta \in End(A, +)$ is *diagonal* when, for each element $e_i$ of the standard basis, there exists $d_i \in \mathbb{Z}_{|e_i|}$ such that $\delta(e_i) = d_i e_i$. If $\varepsilon \in End(A, +)$, a non-zero element $a \in A$ is called a *characteristic element* of $\varepsilon$ if there exists a $\lambda \in \mathbb{Z}_{|a|}$ such that $\varepsilon(a) = \lambda a$. In this case we call $\lambda$ the *characteristic value* of $\varepsilon$ associated with $a$. We say an endomorphism $\varepsilon \in End(A, +)$ is *diagonalizable* if $A$ has a ppc-basis consisting of characteristic elements of $\varepsilon$. The following result is analogous to a familiar one in linear algebra.

**Proposition 5.1:** If $A$ is a finite abelian group and $\varepsilon \in End(A, +)$ is diagonalizable, then $\varepsilon$ is similar to a diagonal endomorphism.

**Proof:** Let us suppose that $\bar{x} = \{x_1, x_2, \ldots, x_r\}$ is a ppc-basis of $A$ with each $x_i$ a characteristic element of $\varepsilon$. Thus for each $x_i$ there is a $\lambda_i \in \mathbb{Z}_{|x_i|}$ such that $\varepsilon(x_i) = \lambda_i x_i$. First, we wish to define an automorphism $\alpha \in Aut(A, +)$ sending the standard basis $\bar{e}$ to this basis $\bar{x}$. Since $\bar{e}$ and $\bar{x}$ are both ppc-bases for $A$, then $A = \langle e_1 \rangle \oplus \cdots \oplus \langle e_r \rangle = \langle x_1 \rangle \oplus \cdots \oplus \langle x_r \rangle$ give decompositions of $A$ into direct sums of cyclic groups of prime power order and, as mentioned earlier, there is a bijection of these summands pairing each $\langle e_i \rangle$ with a unique isomorphic $\langle x_j \rangle$. Thus there is a permutation $\pi$ of $\{1, 2, \ldots, r\}$ so that $\langle e_i \rangle \cong \langle x_{\pi(i)} \rangle$, for each $i$. In particular, $|x_{\pi(i)}|$ equals, and hence divides, $|e_i|$. Thus the mapping $\alpha^* : \bar{e} \to \bar{x}$ defined by $\alpha^*(e_i) = x_{\pi(i)}$ can be extended to an endomorphism, $\alpha$, of $A$. And since it is a bijection of the two ppc-bases of $A$, $\alpha$ is an automorphism of $A$. Next we define a mapping $\delta^* : \bar{e} \to A$ by $\delta^*(e_i) = \lambda_{\pi(i)} e_i$. Since $\lambda_{\pi(i)} \in \mathbb{Z}_{|x_{\pi(i)}|}$, and $|x_{\pi(i)}| = |e_i|$, it follows that $\lambda_{\pi(i)} \in \mathbb{Z}_{|e_i|}$ and, hence, that the extension of this map, $\delta$, is a diagonal endomorphism of $A$. Note now that $\varepsilon\alpha(e_i) = \varepsilon(x_{\pi(i)}) = \lambda_{\pi(i)} x_{\pi(i)}$, while $\alpha\delta(e_i) = \alpha(\lambda_{\pi(i)} e_i)$ $= \lambda_{\pi(i)} \alpha(e_i) = \lambda_{\pi(i)} x_{\pi(i)}$. Thus $\alpha\delta = \varepsilon\alpha$ and it follows that $\alpha^{-1}\varepsilon\alpha = \delta$, and thus, that $\varepsilon$ is similar to $\delta$. ∎

**Lemma 5.2:** If $(A, +)$ is a finite abelian group, $(\varepsilon_1, \varepsilon_2)$ is a pair of commuting endomorphisms of $A$, and $(\eta_1, \eta_2)$ is a pair of endomorphisms of $A$ which is similar to $(\varepsilon_1, \varepsilon_2)$, then $(\eta_1, \eta_2)$ is a commuting pair.

**Proof:** Since $(\varepsilon_1, \varepsilon_2) \sim (\eta_1, \eta_2)$, there is an $\alpha \in Aut(A, +)$ for which $\eta_1 = \alpha^{-1}\varepsilon_1\alpha$ and $\eta_2 = \alpha^{-1}\varepsilon_2\alpha$. Therefore,
$$\eta_1\eta_2 = \alpha^{-1}\varepsilon_1\alpha\alpha^{-1}\varepsilon_2\alpha = \alpha^{-1}\varepsilon_1\varepsilon_2\alpha = \alpha^{-1}\varepsilon_2\varepsilon_1\alpha = \alpha^{-1}\varepsilon_2\alpha\alpha^{-1}\varepsilon_1\alpha = \eta_2\eta_1. \blacksquare$$

**Lemma 5.3:** If $A$ is a finite abelian group, $\varepsilon \in End(A, +)$ is idempotent, and $0 \neq a \in A$ is a prime power order characteristic element of $\varepsilon$ with characteristic value $\lambda$, then $\lambda \in \{0, 1\} \subseteq \mathbb{Z}_{|a|}$.

**Proof:** Suppose $p$ is a prime and $n$ is a positive integer such that $|a| = p^n$. Since there is a $\lambda \in \mathbb{Z}_{p^n}$ such that $\varepsilon(a) = \lambda a$, we have,
$$\lambda a = \varepsilon(a) = \varepsilon^2(a) = \varepsilon(\varepsilon(a)) = \varepsilon(\lambda a) = \lambda\varepsilon(a) = \lambda^2 a.$$
It follows that $(\lambda^2 - \lambda)a = 0$ and, therefore, that $\lambda^2 \equiv \lambda \pmod{p^n}$. Now suppose that $\lambda \neq 0$ and write $\lambda = p^m s$ where $0 \leq m < n$ and $s$ is coprime to $p$. We may rewrite $\lambda^2 \equiv \lambda \pmod{p^n}$ as $p^{2m} s^2 \equiv p^m s \pmod{p^n}$.

Therefore $p^m s \equiv 1 \pmod{p^{n-m}}$ and, since $n - m \geq 1$, it must be that $m = 0$ and $\lambda \equiv s \equiv 1 \pmod{p^n}$. The result follows. ∎

We will need to consider the structure of idempotent endomorphisms of $A$ more carefully. Let $A$ be a finite abelian group, $\varepsilon \in End(A,+)$, and, for each $i$ $(1 \leq i \leq r)$, $A_i \leq A$ such that $A = A_1 \oplus A_2 \oplus \cdots \oplus A_r$. If $\alpha(A_i) \leq A_i$, for each $i$, we say each $A_i$ is *invariant under* $\varepsilon$. For each $i$, we define $\varepsilon_i \in End(A_i, +)$ to be the restriction of $\varepsilon$ to $A_i$; thus for each $a_i \in A_i$, we have $\varepsilon_i(a_i) = \varepsilon(a_i)$, Note that, since $A_i$ is invariant under $\varepsilon$, $\varepsilon_i$ maps into $A_i$. We call $\varepsilon$ *the direct sum of the endomorphisms* $\varepsilon_1, \ldots, \varepsilon_r$, and write $\varepsilon = \varepsilon_1 + \ldots + \varepsilon_r$. Note that $\varepsilon_i$ is not an endomorphism of $A$. For each $a \in A$, there are unique elements $a_i \in A_i$ so that $a = a_1 + \cdots + a_r$ and, according to our definitions,
$$\varepsilon(a) = \varepsilon(a_1) + \cdots + \varepsilon(a_r) = \varepsilon_1(a_1) + \cdots + \varepsilon_r(a_r).$$

**Lemma 5.4:** If $A$ is a finite abelian group, $\varepsilon \in End(A,+)$, and $A$ can be decomposed as $A = A_1 \oplus A_2 \oplus \cdots \oplus A_r$ with each $A_i$ invariant under $\varepsilon$, and $\varepsilon$ is written as the direct sum $\varepsilon = \varepsilon_1 + \ldots + \varepsilon_r$, then $\varepsilon$ is idempotent if and only if, for each $i$ $(1 \leq i \leq r)$, $\varepsilon_i$ is idempotent.

**Proof:** If $a \in A$ and $a = a_1 + \cdots + a_r$, where each $a_i \in A_i$, then $\varepsilon(a) = \varepsilon(a_1) + \cdots + \varepsilon(a_r)$ $= \varepsilon_1(a_1) + \cdots + \varepsilon_r(a_r)$, and $\varepsilon^2(a) = \varepsilon(\varepsilon_1(a_1) + \cdots + \varepsilon_r(a_r)) = \varepsilon_1^2(a_1) + \cdots + \varepsilon_r^2(a_r)$. Therefore, if each $\varepsilon_i$ is idempotent, then $\varepsilon$ is idempotent. Now suppose that some $\varepsilon_i$ were not idempotent. Letting $a' \in A_i$ so that $\varepsilon_i(a') \neq \varepsilon_i^2(a')$, we see that $\varepsilon(a') \neq \varepsilon^2(a')$, and, therefore, $\varepsilon$ is not idempotent. ∎

**Lemma 5.5:** If $A$ is a finite abelian group of ppc-rank $r$, and $\delta$ is a diagonal endomorphism of $A$, then $\delta$ is idempotent if and only if, for each $i$ $(1 \leq i \leq r)$, $\delta(e_i) \in \{0, e_i\}$.

**Proof:** Since $\delta$ is diagonal, there is a $d_i \in \mathbb{Z}_{|e_i|}$, for each $i$ $(1 \leq i \leq r)$, such that $\delta(e_i) = d_i e_i$. Thus each $e_i$ is a characteristic element of $\delta$ and each subgroup $\langle e_i \rangle$ is invariant under $\delta$. Writing $A = \langle e_1 \rangle \oplus \cdots \oplus \langle e_r \rangle$ and decomposing $\delta$ over this sum as $\delta = \delta_1 + \cdots + \delta_r$, Lemma 5.4 implies that if $\delta$ is idempotent then each $\delta_i$ is idempotent. Since each $e_i$ is of prime power order, Lemma 5.3 implies that, $d_i$, the characteristic value of $\delta_i$ associated with $e_i$ is either 0 or 1 modulo $|e_i|$. Thus we have $\delta(e_i) \in \{0, e_i\}$, as required. Conversely, suppose that for each $i$ we have $\delta(e_i) \in \{0, e_i\}$. If $e_i$ has characteristic value 1, then $e_i = \delta(e_i) = \delta_i(e_i)$. It follows that $\delta_i^2(e_i) = \delta_i(e_i) = e_i$ and $\delta_i$ is idempotent. On the other hand, if $e_i$ has characteristic value 0, then $0 = \delta(e_i) = \delta_i(e_i)$ and $\delta_i^2(e_i) = \delta_i(0) = 0$. In either case, we conclude that each $\delta_i$ is idempotent and, by Lemma 5.4, that $\delta$ is therefore also idempotent. ∎

**Lemma 5.6:** If $A$ is a finite abelian group, then idempotent diagonal endomorphisms of $A$ commute.
**Proof:** Let $A$ be of ppc-rank $r$ and let $\delta_1$ and $\delta_2$ be idempotent diagonal endomorphisms of $A$. By Lemma 5.5, for $j = 1, 2$ and for each $i$ $(1 \leq i \leq r)$, we have $\delta_j(e_i) \in \{0, e_i\}$. We claim that $\delta_1 \delta_2(e_i) = \delta_2 \delta_1(e_i)$ for each $i$. There are four cases to consider: (i) $\delta_1(e_i) = \delta_2(e_i) = 0$, (ii) $\delta_1(e_i) = 0, \delta_2(e_i) = e_i$, (iii) $\delta_1(e_i) = e_i, \delta_2(e_i) = 0$, and (iv) $\delta_1(e_i) = \delta_2(e_i) = e_i$. In each case it is routine to check that $\delta_1$ and $\delta_2$ commute when applied to each $e_i$. Writing $a \in A$ as $a = m_1 e_1 + \ldots + m_r e_r$, we see that

$$\delta_1 \delta_2(a) = \sum_{i=1}^{r} m_i \delta_1 \delta_2(e_i) = \sum_{i=1}^{r} m_i \delta_2 \delta_1(e_i) = \delta_1 \delta_2(a). \blacksquare$$

**Lemma 5.7:** If $A$ is a finite abelian group, $\varepsilon$ and $\eta$ are similar endomorphisms of $A$, and $\varepsilon$ is idempotent, then $\eta$ is idempotent.

**Proof:** Let $\alpha \in Aut(A,+)$ so that $\alpha^{-1} \varepsilon \alpha = \eta$. Therefore $\eta^2 = (\alpha^{-1} \varepsilon \alpha)^2 = \alpha^{-1} \varepsilon^2 \alpha = \alpha^{-1} \varepsilon \alpha = \eta$. $\blacksquare$

**Lemma 5.8:** If $(A,+)$ is an abelian group, and $\varepsilon \in End(A,+)$ is idempotent, then $A = \varepsilon(A) \oplus \ker(\varepsilon)$. Furthermore, $\varepsilon(A)$ and $\ker(\varepsilon)$ are invariant under $\varepsilon$ with $\varepsilon$ acting as the identity mapping on $\varepsilon(A)$ and as the zero mapping on $\ker(\varepsilon)$.

**Proof:** For each $a \in A$, $a = \varepsilon(a) + (a - \varepsilon(a))$. Clearly $\varepsilon(a) \in \varepsilon(A)$. To see that $a - \varepsilon(a) \in \ker(\varepsilon)$, note that $\varepsilon(a - \varepsilon(a)) = \varepsilon(a) - \varepsilon^2(a) = 0$. It follows that $A = \varepsilon(A) + \ker(\varepsilon)$. To see that the sum is direct, let $x \in \varepsilon(A) \cap \ker(\varepsilon)$. Since $x \in \varepsilon(A)$, there is an $x' \in A$ with $\varepsilon(x') = x$. Therefore we have,

$$x = \varepsilon(x') = \varepsilon^2(x') = \varepsilon(x).$$

Since $x \in \ker(\varepsilon)$, this shows that $x = 0$ and, therefore, $\varepsilon(A) \cap \ker(\varepsilon) = \{0\}$. Since $\varepsilon$ acts as the identity mapping on its image and the zero mapping on its kernel, it follows that both subgroups are invariant under $\varepsilon$. $\blacksquare$

**Theorem 5.9:** Let $(A,+)$ be a finite abelian group and let $(\varepsilon_1, \varepsilon_2)$ be a pair of endomorphisms of $A$. The pair $(\varepsilon_1, \varepsilon_2)$ is a commuting pair of idempotent endomorphisms of $A$ if and only if there exists a pair $(\delta_1, \delta_2)$ of idempotent diagonal endomorphisms of $A$ similar to $(\varepsilon_1, \varepsilon_2)$.

**Proof:** Suppose first that $(\varepsilon_1, \varepsilon_2)$ is similar to a pair $(\delta_1, \delta_2)$ of idempotent diagonal endomorphisms of $A$. Since $\delta_1$ and $\delta_2$ are idempotent and similar to $\varepsilon_1$ and $\varepsilon_2$, respectively, Lemma 5.7 implies that $\varepsilon_1$ and $\varepsilon_2$ are idempotent. Lemma 5.6 implies that the pair $(\delta_1, \delta_2)$ commutes and, thus, by Lemma 5.2, the pair $(\varepsilon_1, \varepsilon_2)$ commutes.

Conversely, suppose that $(\varepsilon_1, \varepsilon_2)$ is a commuting pair of idempotent endomorphisms of $A$. We will apply Lemma 5.8 to decompose $A$ as $A = \varepsilon_1(A) \oplus \ker(\varepsilon_1)$. Since $\varepsilon_1$ acts as the identity map on $\varepsilon_1(A)$ and the zero mapping on $\ker(\varepsilon_1)$, each non-zero element, $x$, of these summands is a characteristic element of $\varepsilon_1$ having characteristic value 1 if $x \in \varepsilon_1(A)$ and 0 if $x \in \ker(\varepsilon_1)$.

*Claim:* Both $\varepsilon_1(A)$ and $\ker(\varepsilon_1)$ are invariant under $\varepsilon_2$.

*Proof of claim:* Since $\varepsilon_1$ acts as the identity mapping on $\varepsilon_1(A)$, we have $x = \varepsilon_1(x)$, for each $x \in \varepsilon_1(A)$. Applying $\varepsilon_2$ and using the commutativity of the two endomorphisms we have, $\varepsilon_2(x) = \varepsilon_2 \varepsilon_1(x) = \varepsilon_1 \varepsilon_2(x)$. Thus $\varepsilon_2(x) \in \varepsilon_1(A)$ and it follows that $\varepsilon_2$ maps $\varepsilon_1(A)$ into itself. Thus $\varepsilon_1(A)$ is invariant under $\varepsilon_2$. Now let $x \in \ker(\varepsilon_1)$. Since $\varepsilon_1(x) = 0$, it follows that $0 = \varepsilon_2(0) = \varepsilon_2 \varepsilon_1(x) = \varepsilon_1 \varepsilon_2(x)$. Hence $\varepsilon_2(x) \in \ker(\varepsilon_1)$ and this shows that $\ker(\varepsilon_1)$ is invariant under $\varepsilon_2$, and the claim is established.

Thus we may decompose the endomorphism $\varepsilon_2$ according to the direct sum decomposition $A = \varepsilon_1(A) \oplus \ker(\varepsilon_1)$ as $\varepsilon_2 = \varepsilon_{2,1} + \varepsilon_{2,2}$, where $\varepsilon_{2,1}$ denotes the restriction of $\varepsilon_2$ to $\varepsilon_1(A)$ and $\varepsilon_{2,2}$ denotes the restriction of $\varepsilon_2$ to $\ker(\varepsilon_1)$. Since $\varepsilon_2$ is idempotent, Lemma 5.4 implies that $\varepsilon_{2,1}$ and $\varepsilon_{2,2}$ are idempotent. Applying Lemma 5.8 to the endomorphisms $\varepsilon_{2,1}$ and $\varepsilon_{2,2}$ acting on $\varepsilon_1(A)$ and $\ker(\varepsilon_1)$, respectively, we obtain decompositions of these as

$$\varepsilon_1(A) = \varepsilon_{2,1}(\varepsilon_1(A)) \oplus \ker(\varepsilon_{2,1}) \text{ and } \ker(\varepsilon_1) = \varepsilon_{2,2}(\ker \varepsilon_1) \oplus \ker \varepsilon_{2,2}.$$

We observe that the non-zero elements of $\varepsilon_{2,1}(\varepsilon_1(A))$ and $\ker(\varepsilon_{2,1})$ are characteristic elements of $\varepsilon_2$ with $\lambda = 1$ and $\lambda = 0$, respectively. Also the non-zero elements of $\varepsilon_{2,2}(\ker(\varepsilon_1))$ and $\ker(\varepsilon_{2,2})$ are characteristic elements of $\varepsilon_2$ with $\lambda = 1$ and $\lambda = 0$, respectively. Writing

$$A = \varepsilon_{2,1}(\varepsilon_1(A)) \oplus \ker(\varepsilon_{2,1}) \oplus \varepsilon_{2,2}(\ker(\varepsilon_1)) \oplus \ker(\varepsilon_{2,2}), \tag{5.1}$$

we see that the non-zero elements of these four summands are simultaneously characteristic elements of both $\varepsilon_1$ and $\varepsilon_2$. Since each of the four summands is a finite abelian group, each can be given a ppc-decomposition. Let us suppose that the four summands in (5.1) have ppc-bases $\overline{x}_1, \overline{x}_2, \overline{x}_3$, and, $\overline{x}_4$, of ranks $r_1, r_2, r_3$, and $r_4$, respectively with $r = r_1 + r_2 + r_3 + r_4$. If we let $s_1 = r_1, s_2 = s_1 + r_2, s_3 = s_2 + r_3$, and $s_4 = s_3 + r_4 (= r)$, we can number these basis elements consecutively as

$$\overline{x}_1 = \{x_1, x_2, \ldots, x_{s_1}\}, \overline{x}_2 = \{x_{s_1+1}, x_{s_1+2}, \ldots, x_{s_2}\}, \overline{x}_3 = \{x_{s_2+1}, x_{s_2+2}, \ldots, x_{s_3}\}, \overline{x}_4 = \{x_{s_3+1}, x_{s_3+2}, \ldots, x_{s_4}\}.$$

Consider the standard basis, $\overline{e} = \{e_1, \ldots, e_r\}$ of $A$. By the fundamental structure theorem, and as in the proof of Proposition 5.1, we know that there is a permutation $\pi$ of $\{1, 2, \ldots, r\}$ such that $|e_i| = |x_{\pi(i)}|$, for each $i (1 \leq i \leq r)$. It follows that the mapping $\alpha : A \to A$ defined by extending $\alpha(e_i) = x_{\pi(i)}$, for each $i$, is a automorphism of $A$. Since each $x_i$ is a characteristic vector of both $\varepsilon_1$ and $\varepsilon_2$, and since both are idempotent, by Lemma 5.3, there exist $\lambda_i, \mu_i \in \{0, 1\}$ such that $\varepsilon_1(x_i) = \lambda_i x_i$ and $\varepsilon_2(x_i) = \mu_i x_i$. We define two diagonal endomorphisms as follows: $\delta_1(e_i) = \lambda_{\pi(i)} e_i$ and $\delta_2(e_i) = \mu_{\pi(i)} e_i$. It follows then that,

$$\varepsilon_1 \alpha(e_i) = \varepsilon_1 x_{\pi(i)} = \lambda_{\pi(i)} x_{\pi(i)}, \ \alpha \delta_1(e_i) = \alpha(\lambda_{\pi(i)} e_i) = \lambda_{\pi(i)} \alpha(e_i) = \lambda_{\pi(i)} x_{\pi(i)},$$

$$\varepsilon_2 \alpha(e_i) = \varepsilon_2 x_{\pi(i)} = \mu_{\pi(i)} x_{\pi(i)}, \text{ and } \alpha \delta_2(e_i) = \alpha(\mu_{\pi(i)} e_i) = \mu_{\pi(i)} \alpha(e_i) = \mu_{\pi(i)} x_{\pi(i)}.$$

Thus we have $\varepsilon_1 \alpha = \alpha \delta_1$, $\varepsilon_2 \alpha = \alpha \delta_2$; hence, $\alpha^{-1} \varepsilon_1 \alpha = \delta_1$ and $\alpha^{-1} \varepsilon_2 \alpha = \delta_2$, as required. ∎

Note that the pair $(\delta_1, \delta_2)$ constructed in this proof commutes by Lemma 5.6.

**Theorem 5.10:** If $A$ is a finite abelian group of ppc-rank $r$, then there are at most $4^r$ isomorphism classes of associative interchange rings based on $A$. This bound is best possible in the sense that, for each $r > 0$, there exists a finite abelian group, $A_r$, of ppc-rank $r$ with exactly $4^r$ isomorphism classes of associative interchange rings based on $A_r$.

**Proof:** If $A$ is any finite abelian group, Theorem 4.3 implies that each pair of commuting idempotent endomorphisms of $A$ yields an associative interchange ring and all associative interchange rings based on $A$ arise in this manner. Applying Theorem 5.9, we know that each isomorphism class of these interchange rings contains an interchange ring constructed from a pair $(\delta_1, \delta_2)$ of idempotent diagonal endomorphisms of $A$. If $\delta$ is an idempotent diagonal endomorphism of $A$, Lemma 5.5 implies that, for each $i$ $(1 \leq i \leq r)$, $\delta(e_i) \in \{0, e_i\}$. There are $2^r$ such $\delta$'s and, thus, there are $(2^r)^2$ pairs of these. Many of these pairs may be similar and hence the number of isomorphism classes of associative interchange rings based on $A$ may be less, but $4^r$ provides an upper bound.

To show that this bound is best possible in the sense stated, let $p_i$ denote the $i^{th}$ prime number, let $n(r) = p_1 \cdots p_r$, and let $A_r$ be the group $(\mathbb{Z}_{n(r)}, +)$. Note that $\mathbb{Z}_{n(r)} = \mathbb{Z}_{p_1} \oplus \cdots \oplus \mathbb{Z}_{p_r}$ is of ppc-rank $r$ and is generated additively by the multiplicative identity 1 of the ring $(\mathbb{Z}_{n(r)}, +, \cdot)$. If $\varepsilon$ is an endomorphism of $\mathbb{Z}_{n(r)}$ and $x \in \mathbb{Z}_{n(r)}$, then $\varepsilon(x) = x\varepsilon(1)$. Thus for each endomorphism there is an element $a \in \mathbb{Z}_{n(r)}$ so that $\varepsilon(x) = ax$, and it follows that, for each $i$, $\varepsilon(e_i) = ae_i$ when $a$ is reduced modulo $|e_i| = p_i$. Therefore each endomorphism, $\varepsilon$, is diagonal and, by Lemma 5.5, $\varepsilon$ is idempotent if and only if, for each $i$, $\varepsilon(e_i) \in \{0, e_i\}$. Since $\varepsilon$ is of the form $\varepsilon(x) = ax$ for some $a \in \mathbb{Z}_{n(r)}$, we see that $ae_i \in \{0, e_i\}$. Thus, for each $i$, either $a \equiv 0 \pmod{p_i}$ or $a \equiv 1 \pmod{p_i}$. Writing $a = a_1 + \cdots + a_r$ and $\varepsilon = \varepsilon_1 + \cdots + \varepsilon_r$, we then have each $a_i \in \{0, 1\} \subseteq \mathbb{Z}_{p_i}$. Selecting these in all possible ways, we see there are $2^r$ elements $a \in \mathbb{Z}_{n(r)}$ produced. Thus there are $(2^r)^2$ distinct pairs of these endomorphisms. Note also, if $\varepsilon$ and $\eta$ are endomorphisms of $\mathbb{Z}_{n(r)}$ with $\varepsilon(x) = ax$ and $\eta(x) = bx$, then $\varepsilon\eta(x) = abx = bax = \eta\varepsilon(x)$; thus composition of mappings is commutative in $End(\mathbb{Z}_{n(r)}, +)$. It follows that if $\alpha \in Aut(\mathbb{Z}_{n(r)}, +)$ and $\varepsilon \in End(\mathbb{Z}_{n(r)}, +)$, we have $\alpha^{-1}\varepsilon\alpha = \varepsilon$. Therefore similarity classes of endomorphisms of $\mathbb{Z}_{n(r)}$ are singletons. Thus the $4^r$ pairs of endomorphisms we have generated are non-similar and, by Theorem 3.5, each such pair generates an associative interchange near ring based on $A_r$ in a distinct isomorphism class. ∎

**Theorem 5.11:** If $(A, +)$ is a finite abelian group of ppc-rank $r$, then there are at most $2^r$ isomorphism classes of interchange rings based on $A$ whose multiplicative structure forms a band (*i.e.* an idempotent semigroup.) This bound is best possible in the sense that, for each $r > 0$, there exists a finite abelian group, $A_r$, of ppc-rank $r$ with exactly $2^r$ isomorphism classes of associative interchange rings based on $A_r$.

**Proof:** Since a band is a semigroup, we need only consider associative interchange rings. We know that an associative interchange ring based on $A$ is constructed from a pair of commuting, idempotent endomorphisms of $A$. As seen above, this is similar to a pair $(\delta_1, \delta_2)$ of idempotent diagonal endomorphisms of $A$ for which each standard basis element is sent to either itself or zero. For the multiplicative structure to be idempotent, Proposition 4.2 tells us that $\alpha + \beta = \iota$, the identity mapping. Thus we know that $\delta_2 = \delta_1 - \iota$, that is, $\delta_2(e_i) = e_i$ if $\delta_1(e_i) = 0$ and $\delta_2(e_i) = 0$ if $\delta_1(e_i) = e_i$. Since there are two choices for the image of each $e_i$ under $\delta_1$, there are $2^r$ distinct mappings, $\delta_1$. If the pair is to form an idempotent product, $\delta_2$ is determined uniquely by $\delta_1$. Thus there are $2^r$ such pairs and, hence, at most $2^r$ isomorphism classes of interchange rings whose multiplicative structure is a band.

To see that the bound is best possible in the sense claimed, for each $r > 0$, we consider the group $A_r$ discussed in the proof of Theorem 5.10. We saw that there were exactly $2^r$ idempotent diagonal endomorphisms of $A_r$. Pairing each such $\delta$ with $\iota - \delta$ gives a set of $2^r$ pairs which produce, uniquely, all possible isomorphism classes of interchange rings whose multiplicative structure forms a band. ∎

# 6. Interchange rings based on elementary abelian $p^n$-groups.

In the case of an elementary abelian p-group, $A$, we will calculate an exact number of associative interchange rings based on $A$. This number is a cubic polynomial in the ppc-rank of $A$. It turns out that we can carry out this program for a moderately broader class of groups. Given a prime, $p$, and a positive integer, $n$, we call a finite abelian group, $A$, which is a direct sum of $r$ copies of the cyclic group of order $p^n$ an *elementary abelian $p^n$-group* of ppc-rank $r$. Note that when $n=1$, $A$ is an elementary abelian p-group. If $A$ is an elementary abelian $p^n$-group of ppc-rank $r$, then, when $A$ is written as a direct sum using the standard basis, $A = \langle e_1 \rangle \oplus \cdots \oplus \langle e_r \rangle$, each $e_i$ has order $p^n$. In what follows we will maintain the convention that permutations, as well as all other functions, are composed on the left, thus $(\pi_1 \circ \pi_2)(x) = \pi_1(\pi_2(x))$.

**Definition:** Let $p$ be a fixed prime, let $n$ and $r$ be a fixed positive integers, and let $A$ be an elementary abelian $p^n$-group of ppc-rank $r$. An $\alpha \in Aut(A,+)$ is called a *permutation automorphism* of $A$ whenever there exists a permutation, $\pi$, of the set $\{1,\ldots,r\}$ such that $\alpha(e_i) = e_{\pi(i)}$ for each $i$ $(1 \leq i \leq r)$.

Note that for any permutation $\pi$, of the set $\{1,\ldots,r\}$, a permutation automorphism exists as given in the definition. Since $|e_{\pi(i)}| = p^n = |e_i|$, defining $\alpha$ as the extension of $\alpha(e_i) = e_{\pi(i)}$ always yields a homomorphism. It is a bijection since it acts by permuting the standard basis. The permutation automorphism associated with $\pi$ will be denoted $\alpha_\pi$. Note also that the permutation automorphisms of $A$ form a subgroup of $Aut(A,+)$, since $\alpha_\pi^{-1} = \alpha_{\pi^{-1}}$, and if $\pi_1$ and $\pi_2$ are permutations of $\{1,\ldots,r\}$, $\alpha_{\pi_1} \circ \alpha_{\pi_2} = \alpha_{\pi_1 \circ \pi_2}$.

**Lemma 6.1:** If $A$ is an elementary abelian $p^n$-group of ppc-rank $r$, $\alpha_\pi \in Aut(A,+)$ is a *permutation automorphism* of $A$, and $\delta, \eta \in End(A,+)$ with $\delta$ a diagonal endomorphism of $A$ and $\eta$ similar to $\delta$ under $\alpha_\pi$, then $\eta$ is a diagonal endomorphism of $A$.

**Proof:** For each $i$ $(1 \leq i \leq r)$, suppose that $\delta(e_i) = d_i e_i$. Since $\eta = \alpha_\alpha^{-1} \delta \alpha_\alpha$, it follows that
$$\eta(e_i) = \alpha_\pi^{-1} \delta \alpha_\pi(e_i) = \alpha_\pi^{-1} \delta(e_{\pi(i)}) = \alpha_\pi^{-1}(d_{\pi(i)} e_{\pi(i)}) = d_{\pi(i)} \alpha_\pi^{-1}(e_{\pi(i)}) = d_{\pi(i)} e_{\pi^{-1}\pi(i)} = d_{\pi(i)} e_i.$$
Therefore $\eta$ is a diagonal endomorphism of $A$. ∎

If $A$ is an elementary abelian $p^n$-group of ppc-rank $r$, given an integer $s$ $(0 \leq s \leq r)$, we say an idempotent diagonal endomorphism $\delta \in End(A,+)$ is in *s-canonical form* if $\delta(e_i) = e_i$ for $0 \leq i \leq s$ and $\delta(e_i) = 0$ for $s < i \leq r$. Note that if $\delta$ is of this form and $s = 0$, then $\delta$ is the zero map, and if $s = r$, $\delta$ is the identity map. If $\delta$ is an idempotent diagonal endomorphism of an elementary abelian $p^n$-group of ppc-rank $r$, we let $I(\delta) = \{i \in \{1,\ldots,r\} : \delta(e_i) = e_i\}$. Note that $\delta$ is in s-canonical form if and only if $I(\delta) = \{1,\ldots,s\}$.

**Proposition 6.2:** If $A$ is an elementary abelian $p^n$-group of ppc-rank $r$ and $\delta$ is an idempotent diagonal endomorphism of $A$, then there is an idempotent diagonal endomorphism $\eta$ of $A$ in s-canonical form, for $s = |I(\delta)|$, and a permutation automorphism $\alpha_\pi$ under which $\delta$ is similar to $\eta$.

**Proof:** Since $\delta$ is diagonal, for each $i$ $(1 \le i \le r)$, there is a $d_i \in \mathbb{Z}_{|e_i|}$ with $\delta(e_i) = d_i e_i$, and since $\delta$ is idempotent, Lemma 5.3 implies that each $d_i \in \{0,1\}$. If $\delta$ is the zero mapping or the identity mapping, then it is in $s$-canonical form for $s = 0$ and $s = r$, respectively. Selecting $\eta = \delta$ and $\alpha_\pi$ as the identity automorphism, our result follows. Thus we may assume that $I(\delta)$ is a proper, non-empty subset of $\{1, \ldots, r\}$. Let $\pi$ be any permutation of $\{1, \ldots, r\}$ such that $\pi(\{1, \ldots, s\}) = I(\delta)$. We will establish our result by letting $\eta = \alpha_\pi^{-1} \delta \alpha_\pi$ and showing that $\eta$ is in $s$-canonical form. By Lemma 6.1 we know $\eta$ is diagonal. For any $e_i$ we have $\eta(e_i) = \alpha_\pi^{-1} \delta \alpha_\pi(e_i) = \alpha_\pi^{-1} \delta(e_{\pi(i)})$. It is clear that $\pi(i) \in I(\delta)$ if and only if $1 \le i \le s$. Therefore for $1 \le i \le s$,

$$\eta(e_i) = \alpha_\pi^{-1} \delta(e_{\pi(i)}) = \alpha_\pi^{-1}(e_{\pi(i)}) = e_{\pi^{-1}(\pi(i))} = e_i,$$

while for $s+1 \le i \le r$, we have $\pi(i) \notin I(\delta)$ and, since $d_i \in \{0,1\}$, it follows that $\delta(e_{\pi(i)}) = 0$. Therefore,

$$\eta(e_i) = \alpha_\pi^{-1} \delta(e_{\pi(i)}) = \alpha_\pi^{-1}(0) = 0.$$

This establishes that $\eta$ is in $s$-canonical form, and Lemma 5.5 establishes that $\eta$ is idempotent. ∎

Next we will show that, up to similarity, $s$-canonical form is unique.

**Proposition 6.3:** If $A$ is an elementary abelian $p^n$-group, and $\delta$ and $\eta$ are similar idempotent diagonal endomorphisms of $A$ with $\delta$ in $s$-canonical form and $\eta$ in $t$-canonical form, then $s = t$ and, hence, $\eta = \delta$.

**Proof:** Note that if $s = 0$ then $\delta$ is the zero map, which is similar only to itself. Since $\delta \sim \eta$, it follows that $\eta = \delta$ and, therefore, that $t = 0$ and our result follows. Similarly if $t = r$, then $\delta$ is the identity mapping, and the result follows again. Assume then that $1 \le s < t < r$. Since $\delta$ and $\eta$ are similar, there is an $\alpha \in Aut(A, +)$ such that $\alpha \delta = \eta \alpha$. The automorphism $\alpha$ is determined by its action on the standard basis; therefore, for each $i$ there exist $a_{ji} \in \mathbb{Z}_{|e_j|}$ $(1 \le j \le r)$ such that $\alpha(e_i) = \sum_{j=1}^{r} a_{ji} e_j$.

*Claim:* For each $i > s$, $\alpha(e_i) \in \langle e_{t+1}, \ldots, e_r \rangle$.

*Proof of Claim:* Note that for $i > s$, $\alpha \delta(e_i) = \alpha(0) = 0$, while $\eta \alpha(e_i) = \eta\left(\sum_{j=1}^{r} a_{ji} e_j\right) = \sum_{j=1}^{r} a_{ji} \eta(e_j) = \sum_{j=1}^{t} a_{ji} e_j$.

The last equality holds since we have assumed $\eta$ to be in $t$-canonical form. Since $\alpha \delta = \eta \alpha$, this implies, for $i > s$, that $\sum_{j=1}^{t} a_{ji} e_j = 0$. Since the elements $\{e_1, \ldots, e_t\}$ form a ppc-basis for $A$, they form an independent set, and it follows that, for each $i > s$ and $j \le t$, we have $a_{ji} = 0$. Using this fact we see that, for

$i > s$, $\alpha(e_i) = \sum_{j=1}^{r} a_{ji} e_j = \sum_{j=1}^{t} a_{ji} e_j + \sum_{j=t+1}^{r} a_{ji} e_j = \sum_{j=t+1}^{r} a_{ji} e_j$. Thus, for $i > s$, $\alpha(e_i) \in \langle e_{t+1}, \ldots, e_r \rangle$ and our claim is established.

It follows from the claim that the $\alpha \langle e_{s+1}, \ldots, e_r \rangle$ is a subgroup of $\langle e_{t+1}, \ldots, e_r \rangle$. Since $\alpha$ is an automorphism of $A$, it preserves ppc-rank on subgroups; thus we know that

$$r - s = \text{ppc-rank}\langle e_{s+1}, \ldots, e_r \rangle \le \text{ppc-rank}\langle e_{t+1}, \ldots, e_r \rangle = r - t.$$

It follows that $s \geq t$ providing a contradiction to the original assumption that $s < t$. Therefore we can conclude that $s \geq t$. A symmetrical argument shows that $t \geq s$, and it follows that $s = t$. ∎

**Corollary 6.4:** If A is an elementary abelian $p^n$-group of ppc-rank $r$, then there are exactly $r+1$ distinct isomorphism classes of commutative, associative interchange rings based on $A$.
**Proof:** Supposing that a pair, $(\varepsilon, \eta)$, of commuting, idempotent, endomorphisms of $A$ generates a commutative, associative interchange ring. By Proposition 4.1, we have $\varepsilon = \eta$. When we put $\varepsilon$ into $t$-canonical form, $\eta$ is automatically put into the same form. By Propositions 6.2 and 6.3 we see that this generates a distinct isomorphism class of interchange rings for each value of $t$. Since $t$ ranges from 0 to $r$, our result follows. ∎

If $A$ is an elementary abelian $p^n$-group of ppc-rank $r$, $\delta$ is an idempotent diagonal endomorphism of $A$, and $s, t_1$, and $t_2$ are integers such that $0 \leq t_1 \leq s \leq t_2 \leq r$ we say $\delta$ is in $s, t_1, t_2$-*canonical form* when $\delta(e_i) = e_i$ if $0 \leq i \leq t_1$ and $s < i \leq t_2$ and $\delta(e_i) = 0$ if $t_1 < i \leq s$ and $t_2 < i \leq r$. We introduce one further notation; for any integer $s$ $(1 \leq s \leq r)$, let $I_{s,1}(\delta) = \{i \in \{1, \ldots, s\} : \delta(e_i) = e_i\}$ and $I_{s,2}(\delta) = \{i \in \{s+1, \ldots, r\} : \delta(e_i) = e_i\}$.

**Proposition 6.5:** If $A$ is a finite elementary abelian $p^n$-group of ppc-rank $r$, $\delta$ is an idempotent diagonal endomorphism of $A$, and $s$ $(0 \leq s \leq r)$ is an integer, then there exists an idempotent diagonal endomorphism $\eta$ of $A$ in $s, t_1, t_2$-canonical form for $t_1 = |I_{s,1}(\delta)|$, $t_2 = s + |I_{s,2}(\delta)|$, and a permutation automorphism $\alpha_\pi$ such that $\delta$ is similar to $\eta$ under $\alpha_\pi$ and, for $0 < s < r$, the subgroups $A_s = \langle e_1, \ldots, e_s \rangle$ and $B_s = \langle e_{s+1}, \ldots, e_r \rangle$ are invariant under $\alpha_\pi$.
**Proof:** If $s = 0$ then $I_{s,1}(\delta) = \emptyset$ and $I_{s,2}(\delta) = I(\delta)$, and our result follows immediately from Proposition 6.2. If $s = r$, then $I_{s,1}(\delta) = I(\delta)$ and $I_{s,2}(\delta) = \emptyset$, and the result follows, again, from Proposition 6.2. Thus we may assume that $0 < s < r$. Note that, by Lemma 5.5, since $\delta$ is idempotent, for each $i$ we have $\delta(e_i) \in \{0, e_i\}$; therefore, $A_s$ and $B_s$ are both invariant under $\delta$. Since $A = A_s \oplus B_s$ we can then write $\delta$ as the direct sum $\delta = \delta_1 + \delta_2$, where $\delta_1$ and $\delta_2$ are the restrictions of $\delta$ to $A_s$ and $B_s$ respectively. Note that by Lemma 5.4, $\delta_1$ and $\delta_2$ are idempotent. Since $\delta(e_i) \in \{0, e_i\}$ the same is true of $\delta_1$ and $\delta_2$ thus they are diagonal as well. We then apply Proposition 6.2 to both $\delta_1 \in End(A_s, +)$ and $\delta_2 \in End(B_s, +)$. We then have idempotent diagonal endomorphisms $\eta_1 \in End(A_s, +)$ in $t_1$-canonical form and $\eta_2 \in End(B_s, +)$ in $t_2'$-canonical form, for $t_2' = t_2 - s (= |I_{s,2}(\delta)|)$, and permutation automorphisms $\alpha_{\pi_1} \in Aut(A_s, +)$ and $\alpha_{\pi_2} \in Aut(B_s, +)$ so that $\delta_1$ is similar to $\eta_1$ under $\alpha_{\pi_1}$ and $\delta_2$ is similar to $\eta_2$ under $\alpha_{\pi_2}$.
Let $\eta = \eta_1 + \eta_2 (\in End(A, +))$ and $\alpha = \alpha_{\pi_1} + \alpha_{\pi_2}$ $(\in Aut(A, +))$, and note that since $\alpha$ is defined as a direct product, the subgroups $A_s$ and $B_s$ are invariant under $\alpha$. To complete the proof we need to establish that (i) $\eta$ is in $s, t_1, t_2$-canonical form, (ii) $\delta$ is similar to $\eta$ under $\alpha$, and that (iii) $\alpha$ is a permutation automorphism. To prove (i) we first calculate $\eta(e_i)$ for $1 \leq i \leq s$. Since $e_i \in A_s$, $\eta(e_i) = \eta_1(e_i)$. Since $\eta_1$ is in $t_1$-canonical form, we have $\eta(e_i) = e_i$ for $1 \leq i \leq t_1$ and $\eta(e_i) = 0$ for $t_1 < i \leq s$. Now considering $\eta(e_i)$ for $s < i \leq r$, we have $\eta(e_i) = \eta_2(e_i)$. Since $\eta_2$ is in $t_2'$-canonical form, we obtain $\eta(e_i) = e_i$ for $s < i \leq t_2$ and $\eta(e_i) = 0$ for $t_2 < i \leq r$. Thus $\eta$ is in $s, t_1, t_2$-canonical form and (i) is established. To prove (ii) we note that
$$\alpha^{-1}\delta\alpha = (\alpha_1 + \alpha_2)^{-1}(\delta_1 + \delta_2)(\alpha_1 + \alpha_2) = \alpha_1^{-1}\delta_1\alpha_1 + \alpha_2^{-1}\delta_2\alpha_2 = \eta_1 + \eta_2 = \eta.$$
To establish (iii) let $\pi$ be the permutation of $\{1,\ldots,r\}$ which acts as $\pi_1$ on $\{1,\ldots,s\}$ and as $\pi_2$ on $\{s+1,\ldots,r\}$. Thus we have $\alpha = \alpha_{\pi_1} + \alpha_{\pi_2} = \alpha_\pi$ a permutation automorphism of $A$. ∎

If $(\delta_1, \delta_2)$ is a pair of idempotent diagonal endomorphisms of $A$, $s, t_1,$ and $t_2$ are integers with $1 \leq t_1 \leq s \leq t_2 \leq r$ such that $\delta_1$ is in $s$-canonical form and $\delta_2$ in $s, t_1, t_2$-canonical form, then we say that the pair $(\delta_1, \delta_2)$ is in *canonical form*.

**Proposition 6.6:** If $A$ is an elementary abelian $p^n$-group of ppc-rank $r$, $(\delta_1, \delta_2)$ is a pair of idempotent diagonal endomorphisms of $A$, then for $s = |I(\delta_1)|$, $t_1 = |I_{s,1}(\delta_2)|$, and $t_2 = |I_{s,2}(\delta_2)|$, there exists a pair of idempotent diagonal endomorphisms of $A$, $(\eta_1, \eta_2)$, such that $(\delta_1, \delta_2) \sim (\eta_1, \eta_2)$, $\eta_1$ is in $s$-canonical form, and $\eta_2$ is in $s, t_1, t_2$-canonical form. In other words, any pair of idempotent diagonal endomorphism of $A$ can be put into canonical form.

**Proof:** Focussing first on $\delta_1$, by Proposition 6.2, we know that there is an idempotent diagonal endomorphism $\eta_1$ of $A$ in $s$-canonical form with $\delta_1$ similar to $\eta_1$ under a permutation automorphism $\alpha$. Letting $\delta_2' = \alpha^{-1}\delta_2\alpha$ we see, by Lemma 6.1, that $\delta_2'$ is diagonal and, by Lemma 5.7, $\delta_2'$ is idempotent. Thus we have $(\delta_1, \delta_2)$ similar under a permutation automorphism to an idempotent diagonal pair $(\eta_1, \delta_2')$ with $\eta_1$ in $s$-canonical form. Applying Proposition 6.5 to $\delta_2'$ with $s = |I(\delta_1)|$, we have the existence of an idempotent diagonal endomorphism $\eta_2$ of $A$ in $s, t_1, t_2$-canonical form and a permutation automorphism $\alpha_\pi$ such that $\delta_2'$ is similar to $\eta_2$ under $\alpha_\pi$ with the subgroups $A_s = \langle e_1, \ldots, e_s \rangle$ and $B_s = \langle e_{s+1}, \ldots, e_r \rangle$ invariant under $\alpha_\pi$.

*Claim:* $\alpha_\pi$ commutes with $\eta_1$.

*Proof of Claim:* Note that since $A_s$ is invariant under $\alpha_\pi$, for $1 \leq i \leq s$, we have $e_{\pi(i)} = \alpha_\pi(e_i) \in A_s$. It follows that $\eta_1(e_{\pi(i)}) = e_{\pi(i)}$. Therefore $\eta_1\alpha_\pi(e_i) = \eta_1(e_{\pi(i)}) = e_{\pi(i)}$, while $\alpha_\pi\eta_1(e_i) = \alpha_\pi(e_i) = e_{\pi(i)}$. Thus $\alpha_\pi$ commutes with $\eta_1$ on $A_s$. Since $B_s$ is invariant under $\alpha_\pi$, for $s < i \leq r$, we have $e_{\pi(i)} = \alpha_\pi(e_i) \in B_s$. Thus $\eta_1(e_{\pi(i)}) = 0$. Thus we have $\eta_1\alpha_\pi(e_i) = \eta_1(e_{\pi(i)}) = 0$, while $\alpha_\pi\eta_1(e_i) = \alpha_\pi(0) = 0$. Thus we see that $\alpha_\pi$ commutes with $\eta_1$ on $B_s$ and the claim follows.

Since $(\delta_1, \delta_2) \sim (\eta_1, \delta_2')$, to complete the proof we must show that $(\eta_1, \delta_2')$ is similar to $(\eta_1, \eta_2)$. Since $\alpha_\pi$ commutes with $\eta_1$, we have $\alpha_\pi^{-1}(\eta_1, \delta_2')\alpha_\pi = (\alpha_\pi^{-1}\eta_1\alpha_\pi, \alpha_\pi^{-1}\delta_2'\alpha_\pi) = (\eta_1, \eta_2)$. ∎

**Proposition 6.7:** If $A$ is an elementary abelian $p^n$-group of ppc-rank $r$, $(\delta_1, \delta_2)$ and $(\eta_1, \eta_2)$ are pairs of idempotent diagonal endomorphisms of $A$, $s, t_1, t_2, s', t_1',$ and $t_2'$ are integers with $1 \leq t_1 \leq s \leq t_2 \leq r$ and $1 \leq t_1' \leq s' \leq t_2' \leq r$, with $\delta_1$ in $s$-canonical form, $\delta_2$ in $s, t_1, t_2$-canonical form, $\eta_1$ in $s'$-canonical form, and $\eta_2$ in $s', t_1', t_2'$-canonical form, and $(\delta_1, \delta_2) \sim (\eta_1, \eta_2)$, then $s = s', t_1 = t_1'$, and $t_2 = t_2'$. Thus if both pairs are in canonical form and similar, they are equal.

**Proof:** Let $\alpha \in Aut(A, +)$ so that $\alpha^{-1}(\delta_1, \delta_2)\alpha = (\eta_1, \eta_2)$. Since $\delta_1$ is similar to $\eta_1$, Proposition 6.3 implies that $s = s'$ and $\delta_1 = \eta_1$. Since $\alpha^{-1}\delta_1\alpha = \eta_1 = \delta_1$, it follows that $\alpha$ commutes with $\delta_1$. Henceforth we will write $s'$ as $s$ and $\eta_1$ as $\delta_1$. For each $i$ and $j$ $(1 \leq i, j \leq r)$, let $a_{ji} \in \mathbb{Z}_{|e_i|}$ so that $\alpha(e_i) = \sum_{j=1}^{r} a_{ji} e_j$.

*Claim:* $\alpha$ is invariant on $A_s$ and $B_s$.

*Proof of claim:* Note that for $1 \leq i \leq s$, $\delta_1(e_i) = e_i$. Since $\alpha$ commutes with $\delta_1$, we have

$$\alpha(e_i) = \alpha\delta_1(e_i) = \delta_1\alpha(e_i) = \delta_1\left(\sum_{i=1}^{r} a_{ji} e_j\right) = \sum_{i=1}^{r} a_{ji}\delta_1(e_j) = \sum_{i=1}^{s} a_{ji}\delta_1(e_j) + \sum_{i=s+1}^{r} a_{ji}\delta_1(e_j)$$

$$= \sum_{i=1}^{s} a_{ji} e_j + \sum_{i=s+1}^{r} a_{ji}(0) = \sum_{j=1}^{s} a_{ji} e_j \in A_s.$$

For each $a \in A_s$, there exist $c_i \in \mathbb{Z}_{|e_i|}$ such that $a = \sum_{i=1}^{s} c_i e_i$. Therefore $\alpha(a) = \alpha\left(\sum_{i=1}^{s} c_i e_i\right) = \sum_{i=1}^{s} c_i \alpha(e_i) \in A_s$.

Next suppose that $s < i \leq r$. We then have $\delta_1(e_i) = 0$. Since $\alpha$ commutes with $\delta_1$, we obtain

$$0 = \alpha(0) = \alpha\delta_1(e_i) = \delta_1\alpha(e_i) = \delta_1\left(\sum_{j=1}^{r} a_{ji} e_j\right) = \sum_{j=1}^{r} a_{ji}\delta_1(e_j) = \sum_{j=1}^{s} a_{ji}\delta_1(e_j) + \sum_{j=s+1}^{r} a_{ji}\delta_1(e_j)$$

$$= \sum_{j=1}^{s} a_{ji} e_j + \sum_{j=s+1}^{r} a_{ji}(0) = \sum_{j=1}^{s} a_{ji} e_j.$$

Since $\sum_{j=1}^{s} a_{ji} e_j = 0$, we conclude that for $s < i \leq r$ and $1 \leq j \leq s$, we have $a_{ji} = 0$. Thus for $s < i \leq r$,

$$\alpha(e_i) = \sum_{j=1}^{s} a_{ji} e_j + \sum_{j=s+1}^{r} a_{ji} e_j = \sum_{j=1}^{s} (0) e_j + \sum_{j=s+1}^{r} a_{ji} e_j = \sum_{j=s+1}^{r} a_{ji} e_j \in B_s.$$

Thus the claim is established.

Note that any diagonal endomorphism is invariant on both $A_s$ and $B_s$ since it sends $e_i$ to some $d_i e_i$. Thus we may decompose $\alpha, \delta_2,$ and $\eta_2$ over the direct sum $A = A_s \oplus B_s$ as $\alpha = \alpha_1 + \alpha_2$, $\delta_2 = \delta_{2,1} + \delta_{2,2}$, and $\eta_2 = \eta_{2,1} + \eta_{2,2}$. Since $\alpha^{-1}\delta_2\alpha = \eta_2$, we can rewrite this as $(\alpha_1 + \alpha_2)^{-1}(\delta_{2,1} + \delta_{2,2})(\alpha_1 + \alpha_2) = \eta_{2,1} + \eta_{2,2}$. Thus we conclude that $\alpha_1^{-1}\delta_{2,1}\alpha_1 = \eta_{2,1}$ and $\alpha_1^{-1}\delta_{2,2}\alpha_1 = \eta_{2,2}$. Note that $\delta_{2,1}$ and $\eta_{2,1}$ are idempotent diagonal endomorphisms of $A_s$, an elementary abelian $p^n$-group of ppc-rank $s$, which are similar under $\alpha_1 \in Aut(A_s, +)$. Note also that $\delta_{2,1}$ is in $t_1$-canonical form and $\eta_{2,1}$ is in $t_1'$-canonical form. Applying Proposition 6.3 to we see that $t_1 = t_1'$ and, in fact, $\delta_{2,1} = \eta_{2,1}$. Applying a similar argument to $\delta_{2,2}$ and $\eta_{2,2}$, which are idempotent diagonal endomorphisms of $B_s$ in $t_2$- and $t_2'$-canonical form, respectively, we see that $t_2 = t_2'$ and $\delta_{2,2} = \eta_{2,2}$. It follows that $\delta_2 = \delta_{2,1} + \delta_{2,2} = \eta_{2,1} + \eta_{2,2} = \eta_2$ as required.

**Corollary 6.8:** If $(A,+)$ is an elementary abelian $p^n$-group of ppc-rank $r$ and $(\varepsilon_1, \varepsilon_2)$ is a pair of commuting, idempotent endomorphisms of $A$, then there exists a unique pair $(\delta_1, \delta_2)$ of commuting, idempotent diagonal endomorphisms which is similar to $(\varepsilon_1, \varepsilon_2)$ and in canonical form.

**Proof:** Lemmas 5.9 and Proposition 6.5 show that $(\varepsilon_1, \varepsilon_2)$ is similar to a pair $(\delta_1, \delta_2)$ of idempotent diagonal endomorphisms in canonical form. Lemma 5.6 shows that $\delta_1$ and $\delta_2$ commute. Proposition 6.6 shows that the pair is unique. ∎

**Theorem 6.9:** If $(A,+)$ is an elementary abelian $p^n$-group of ppc-rank $r$, then there are exactly $\frac{1}{6}(r+1)(r+2)(r+3)$ isomorphism classes of associative interchange rings based on $A$.

**Proof:** Theorem 4.3 says that each associative, interchange ring based on $A$ is formed from a pair, $(\varepsilon_1, \varepsilon_2)$ of commuting idempotent endomorphisms of $A$. These pairs may yield isomorphic interchange rings; however, Theorem 3.4 shows that isomorphism of these interchange rings is determined exactly up to similarity of these pairs. Corollary 6.8 implies that there is exactly one pair, $(\delta_1, \delta_2)$, of diagonal endomorphisms in canonical form which is similar to $(\varepsilon_1, \varepsilon_2)$. Thus the isomorphism class of the associative interchange ring generated from $(\varepsilon_1, \varepsilon_2)$ is the same as that for $(\delta_1, \delta_2)$, and $(\delta_1, \delta_2)$ is the only pair in canonical form which generates an element of that isomorphism class. Note that Lemma 5.6 implies that $(\delta_1, \delta_2)$ is a commuting pair. Thus we can count the distinct isomorphism classes of associative interchange rings based on $A$ by counting the number of distinct pairs $(\delta_1, \delta_2)$ of diagonal endomorphism in canonical form. Since $\delta_1$ is in $s$-canonical form for some $s$ ($1 \leq s \leq r$), we see there are $r+1$ possibilities for $\delta_1$. For each of these the $\delta_2$ is in $s, t_1, t_2$-canonical form with $1 \leq t_1 \leq s \leq t_2 \leq r$. With $s$ fixed we can vary $t_1$ from 0 to $s$ and, independently, vary $t_2$ from $s$ to $r$. Thus for each choice of $\delta_1$ there are $(s+1)(r-s+1)$ choices for $\delta_2$. Thus we have exactly $\sum_{s=0}^{r}(s+1)(r-s+1)$ distinct pairs $(\delta_1, \delta_2)$, and this corresponds to the number of distinct isomorphism classes of associative interchange rings based on $A$. Induction on $r$ shows that this bound has the lovely closed form $\frac{1}{6}(r+1)(r+2)(r+3)$. ∎

There are eight isomorphism classes of rings based on the Klein four-group, $V$. In the next example we will find that there are exactly ten isomorphism classes of associative interchange rings based on $V$.

**Example 6.10: Interchange rings based on the Klein four-group.** The Klein four-group is $V = \{0,1,2,3\}$ with Cayley table,

| + | 0 | 1 | 2 | 3 |
|---|---|---|---|---|
| 0 | 0 | 1 | 2 | 3 |
| 1 | 1 | 0 | 3 | 2 |
| 2 | 2 | 3 | 0 | 1 |
| 3 | 3 | 2 | 1 | 0 |

We denote the endomorphism $\varepsilon : V \to V$ sending $0 \mapsto 0, 1 \mapsto x, 2 \mapsto y$, and $3 \mapsto z$ by $\varepsilon = (0xyz)$. $End(V,+)$ has order sixteen with six automorphisms and ten proper endomorphisms.

$End(V,+) = \{\alpha_0 = (0123), \alpha_1 = (0132), \alpha_2 = (0213), \alpha_3 = (0231), \alpha_4 = (0312), \alpha_5 = (0321), \varepsilon_0 = (0000),$
$\varepsilon_1 = (0011), \varepsilon_2 = (0022), \varepsilon_3 = (0033), \varepsilon_4 = (0101), \varepsilon_5 = (0110), \varepsilon_6 = (0202), \varepsilon_7 = (0220),$
$\varepsilon_8 = (0303), \varepsilon_9 = (0330)\}.$

Since $(V,+)$ is abelian, each pair of endomorphisms is image-commuting. The set of idempotents in $End(V,+)$ is $E = \{\alpha_0, \varepsilon_0, \varepsilon_2, \varepsilon_3, \varepsilon_4, \varepsilon_5, \varepsilon_7, \varepsilon_8\}$. We wish to list a complete set of non-similar pairs of commuting idempotents to exhaust the isomorphism classes of associative interchange rings based on $V$. We may begin with $\varepsilon_0$ and $\alpha_0$, which commute with every endomorphism and each have a singleton similarity class. Thus we will include the commuting pairs $(\varepsilon_0, \varepsilon_0), (\alpha_0, \alpha_0), (\alpha_0, \varepsilon_0), (\varepsilon_0, \alpha_0)$ in our list. A straightforward calculation shows that the endomorphisms $E' = \{\varepsilon_2, \varepsilon_3, \varepsilon_4, \varepsilon_5, \varepsilon_7, \varepsilon_8\}$ are similar to each other. Thus, under similarity, the pairs of the form $(\gamma_0, \varepsilon_i)$ with $\gamma_0 \in \{\alpha_0, \varepsilon_0\}$, reduce to $(\alpha_0, \varepsilon_2), (\varepsilon_2, \alpha_0), (\varepsilon_0, \varepsilon_2), (\varepsilon_2, \varepsilon_0)$. Thus we will include these four pairs as well.

We find that $\varepsilon_2\varepsilon_4 = \varepsilon_4\varepsilon_2$, $\varepsilon_3\varepsilon_5 = \varepsilon_5\varepsilon_3$, and $\varepsilon_7\varepsilon_8 = \varepsilon_8\varepsilon_7$, are the only pairs from $E'$ which commute. Furthermore we have $\alpha_1^{-1}(\varepsilon_2, \varepsilon_4)\alpha_1 = (\varepsilon_3, \varepsilon_5)$ and $\alpha_5^{-1}(\varepsilon_2, \varepsilon_4)\alpha_5 = (\varepsilon_7, \varepsilon_8)$. Thus the only candidates for non-similar pairs here are $(\varepsilon_2, \varepsilon_4)$ and $(\varepsilon_4, \varepsilon_2)$. Thus we have found a complete set of ten pairs of commuting idempotent endomorphisms of $(V,+)$ which are non-similar and represent each isomorphism class of associative interchange rings based on the Klein four-group. Among the ten semigroups produced, $(\varepsilon_0, \varepsilon_0)$ generates the zero semigroup, $(\alpha_0, \alpha_0)$ produces an improper interchange ring with the product identical to the sum, and $(\alpha_0, \varepsilon_0)$ and $(\varepsilon_0, \alpha_0)$ generate the left- and right-zero semigroups. There are the four essential interchange rings guaranteed by Theorem 3.6. The remaining six pairs produce more interesting products. For example selection of the pair $(\varepsilon_7, \varepsilon_7)$ yields the interchange ring $(V,+,\bullet)$ with multiplication having the Cayley table below.

| $\bullet$ | 0 | 1 | 2 | 3 |
|---|---|---|---|---|
| 0 | 0 | 2 | 2 | 0 |
| 1 | 2 | 0 | 0 | 2 |
| 2 | 2 | 0 | 0 | 2 |
| 3 | 0 | 2 | 2 | 0 |

Note that for $p = 2$ and $n = 1$, $V$ is a $p^n$-group of ppc-rank $r = 2$; thus Theorem 6.9 predicts exactly $\frac{1}{6}(2+1)(2+2)(3+2) = 10$ isomorphism classes of associative interchange rings and this is what we have found.

**7. Interchange ring theory.**

Since interchange rings are as yet uninvestigated, it would seem natural to ask how much of standard ring theory would transfer over to interchange rings. In this section we will make some elementary first steps in this direction. For full generality we will allow the additive group to be nonabelian, thus we state our theorems for interchange near rings and they will hold for interchange rings as well. As in universal algebra it is natural to identify sub-objects as follows.

**Definition.** If $(R,+,\bullet)$ is an interchange near ring, and $R'$ is a nonempty subset of $R$, we say $R'$ is an *interchange near subring* of $R$, denoted $R' \leq R$, when $(R',+)$ is a subgroup of $(R,+)$ and $(R',\bullet)$ is a submagma of $(R,\bullet)$.

As with rings, the trivial interchange ring $\{0\}$ and the whole of $R$ are interchange near subrings of $R$. To see a nontrivial example, we may take $(R,+,\bullet)$ to be one of the associative interchange rings found in Example 6.10. Here the additive structure is the Klein four-group and the multiplication is defined as $x \bullet y = \varepsilon_7(x+y)$. Letting $R' = \{0,2\}$ it is easy to see that $R' \leq R$. Note that this interchange subring is improper. To see an example of a proper interchange near subring, we refer to Example 3.7, where the additive structure is $(S_3,+)$ and the product is defined as $x \bullet y = \alpha_3(x) + \varepsilon_0(y) = \alpha_3(x)$. Note that $\{0,1,2\}$ is a proper, non-zero interchange near subring of $(S_3,+)$.

In groups and rings congruences correspond to sub-objects, namely normal subgroups and ideals, while for semigroups and lattices this is not the case. It is interesting to notice that for interchange rings congruences correspond to special sub-objects. Thus we can form quotient interchange rings using these.

**Definition:** If $(R,+,\bullet)$ be an interchange near ring, a subset $I$ of $R$ is an *ideal* of $R$, denoted $I \triangleleft R$, if and only if $(I,+,\bullet)$ is an interchange near subring of $(R,+,\bullet)$ and $(I,+)$ is a normal subgroup of $(R,+)$.

We form a congruence relation, $\mathcal{C}$, on $(R,+)$ using the ideal, $I$, in the standard way: for each $x,y \in I$, $(x,y) \in \mathcal{C}$ if and only if $x - y \in I$.

**Theorem 7.1.** Let $(R,+,\bullet)$ be an interchange near ring, and $I \triangleleft R$, then the relation, $\mathcal{C}$, is a congruence relation on $(R,+,\bullet)$.

**Proof:** Since we have selected $N$ as a normal subgroup of $R$, it follows that $\mathcal{C}$ is an equivalence relation on $R$ and that it acts as a congruence on $(R,+)$. It remains to show that it is a congruence on $(R,\bullet)$. So let $x_1, x_2, y_1, y_2 \in R$ so that $x_1 - x_2, y_1 - y_2 \in I$. Therefore there exist $i, i' \in I$ such that $x_1 = x_2 + i$ and $y_1 = y_2 + i'$. Therefore, by the interchange law, $x_1 \bullet y_1 = (x_2 + i) \bullet (y_2 + i') = (x_2 \bullet y_2) + (i \bullet i')$. By hypothesis $i \bullet i' \in I$, thus we have $x_1 \bullet y_1 - x_2 \bullet y_2 \in I$, as required. ∎

Universal algebraic results imply that we can form a unique quotient interchange near ring, $R/I$. It is routine to check that the fundamental homomorphism theorem and the three isomorphism theorems hold in this case. It is also possible to prove, by a straightforward analogy with ring theory, that $I$ is maximal if and only if $R/I$ is simple.

Further points of departure are the investigation of structures like polynomial rings and group rings using interchange rings. It is interesting to note that when one forms $n \times n$ matrices over an interchange (near) ring, $M_n(R)$, a natural definition of matrix addition and multiplication make this into an interchange (near) ring as well.

**Definition:** If $(R,+,\bullet)$ is an interchange (near) ring and $M_n(R)$ is the set of $n \times n$ matrices over $R$, then we may define the *sum* and *product* of matrices $A = (a_{ij})$ and $B = (b_{ij})$ in $M_n(R)$, respectively, as

$$A + B = (a_{ij} + b_{ij}) \text{ and } A \bullet B = \left( \sum_{k=1}^{n} a_{ik} \bullet b_{kj} \right).$$

Note first that in an interchange (near) ring, it follows by a simple induction on $n$ that

$$x_1 \bullet y_1 + x_2 \bullet y_2 + \ldots + x_n \bullet y_n = (x_1 + x_2 + \ldots + x_n) \bullet (y_1 + y_2 + \ldots + y_n).$$

**Theorem 7.2:** If $(R,+,\bullet)$ is an interchange (near) ring, then $(M_n(R),+,\bullet)$ is an interchange (near) ring.

**Proof:** Clearly $(M_n(R),+)$ is a group and is abelian if $(R,+)$ is abelian. It remains to verify the interchange law.

$$A \bullet B + C \bullet D = \left( \sum_{k=1}^{n} a_{ik} \bullet b_{kj} \right) + \left( \sum_{k=1}^{n} c_{ik} \bullet d_{kj} \right) = \sum_{k=1}^{n} \left( a_{ik} \bullet b_{kj} + c_{ik} \bullet d_{kj} \right)$$

$$= \sum_{k=1}^{n} \left( (a_{ik} + c_{kj}) \bullet (b_{ik} + d_{kj}) \right) = (A+C) \bullet (B+D). \blacksquare$$

Thus there are many avenues open to explore

## 8. References.

Mathematics Department, Mount Saint Vincent University, Halifax, Nova Scotia B3M 2J6

email: cedmunds@eastlink.ca